\documentclass[11pt]{article}
 \usepackage{amsfonts}
 \usepackage{amsmath}

 \newcommand{\beq}{\begin{equation}}
 \newcommand{\eeq}{\end{equation}}
 \newcommand{\bea}{\begin{eqnarray}}
 \newcommand{\eea}{\end{eqnarray}}
\newcommand{\beas}{\begin{eqnarray*}}
\newcommand{\eeas}{\end{eqnarray*}}

\newtheorem{theorem}{Theorem}[section]
\newtheorem{definition}[theorem]{Definition}
\newtheorem{proposition}[theorem]{Proposition}
\newtheorem{corollary}[theorem]{Corollary}
\newtheorem{lemma}[theorem]{Lemma}
\newtheorem{remark}[theorem]{Remark}
\newtheorem{example}[theorem]{Example}
\newtheorem{examples}[theorem]{Examples}
\newtheorem{foo}[theorem]{Remarks}

\newenvironment{Remark}{\begin{remark}\rm}{\end{remark}}
\newenvironment{Remarks}{\begin{foo}\rm}{\end{foo}}
\newenvironment{proof}{\addvspace{\medskipamount}\par\noindent{\it Proof}.}
{\unskip\nobreak\hfill$\Box$\par\addvspace{\medskipamount}}

\newcommand{\ang}[1]{\left<#1\right>}  % angular brackets for projection
\newcommand{\brak}[1]{\left(#1\right)}    % round brackets
\newcommand{\crl}[1]{\left\{#1\right\}}   % curly brackets
\newcommand{\edg}[1]{\left[#1\right]}     % edgy brackets

\newcommand{\E}[1]{{\rm E}\left[#1\right]}

\newcommand{\N}[1]{||#1||}     % Norm
\newcommand{\abs}[1]{\left|#1\right|}     % absolute value

 \title{Dynamic Monetary Risk Measures for Bounded
 Discrete-Time Processes}

 \author{Patrick Cheridito
 \\\small ORFE
 \\\small Princeton University
 \\\small Princeton, NJ 08544
 \\\small USA
 \and Freddy Delbaen\thanks{Supported by Credit Suisse}
 \\\small Departement f\"ur Mathematik
 \\\small ETH Z\"urich
 \\\small 8092 Z\"urich
 \\\small Switzerland
 \and Michael Kuppe$\mbox{r}^{*}$
 \\\small Departement f\"ur Mathematik
 \\\small ETH Z\"urich
 \\\small 8092 Z\"urich
 \\\small Switzerland}

 \date{October 15, 2004}

\begin{document}

 \maketitle

 \begin{abstract}
 \noindent
 We study time-consistency questions for processes of monetary risk
 measures that depend on bounded discrete-time processes
 describing the evolution of financial values. The time horizon
 can be finite or infinite. We call a process of monetary risk
 measures time-consistent if it assigns to a process of financial
 values the same risk irrespective of whether it is calculated
 directly or in two steps backwards in time, and we show how this property manifests
 itself in the corresponding process of acceptance sets.
 For processes of coherent and convex monetary risk measures
 admitting a robust representation with sigma-additive linear functionals,
 we give necessary and sufficient conditions for time-consistency in
 terms of the representing functionals.\\[2mm]
 \end{abstract}

 \section{Introduction}

 The notion of coherent risk measure was introduced in
 Artzner et al. (1997, 1999) and further developed in Delbaen (2001, 2002).
 In F\"ollmer and Schied (2002a, b, c) the more general concepts of monetary
 and convex monetary risk measures were introduced. All these
 works discuss one-period risk measurement, that is, the risky
 objects are real-valued random variables describing future
 financial values and the risk of such financial values is only
 measured at the beginning of the time-period considered.
 Some typical examples of financial values in the context of risk measurement are:

 - the market value of a firm's equity

 - the accounting value of a firm's equity

 - the market value of a portfolio of financial securities

 - the surplus of an insurance company\\
 Cvitani\'{c} and Karatzas study the dynamics of a risk associated
 with hedging a given liability in a continuous-time setup.
 In Artzner et al. (2002) the evolution of financial values over time
 is modelled with discrete-time stochastic processes and
 two special classes of time-consistent processes of coherent
 risk measures related to m-stable sets of probability measures
 are introduced. A treatment of the same two classes of time-consistent
 processes of coherent risk measures in continuous time
 and more on m-stable sets can be found in Delbaen (2004).
 Engwerda et al. (2002) is similar to Artzner et al. (2002)
 but also discusses the effects of hedging and
 the applicability of dynamic programming algorithms.
 Cheridito et al. (2004a, b) contain representation results for
 coherent and convex monetary risk measures that depend on
 processes of financial values evolving in continuous time.
 Rosazza Gianin (2003) studies the
 relation between risk measures and g-expectations.
 Frittelli and Rosazza Gianin (2004) contains a summary of earlier
 results on convex monetary risk measures and connections to
 indifference pricing and g-expectations.
 Riedel (2004), Detlefsen (2003), Scandolo (2003) and Weber (2003)
 study dynamic coherent or convex monetary risk measures for
 cash-flow streams in discrete time.

 In this paper we follow Artzner et al. and measure the risk of
 discrete-time processes of financial values. We simply call them
 value processes. Of course, in
 discrete-time, value processes can be turned into cash-flow
 streams by passing to increment processes. But this
 transformation does not preserve the order of
 almost sure dominance, and because this order plays a crucial role in
 our definition of monetary risk measures, it makes a difference
 whether we take the risky objects to be value processes or
 cash-flow streams. Since in most practical applications it
 can be assumed that money can be borrowed and lent at a risk-free rate,
 we find it more natural to order value processes than cash-flow streams
 by almost sure dominance.

 The structure of the paper is as follows. In Section 2 we
 introduce the basic setup and some notation. In Section 3 we
 introduce monetary risk measures conditioned on the information
 available at stopping times, study the
 relation between such risk measures and their acceptance sets and
 prove conditioned representation results for coherent and convex
 monetary risk measures. In Section 4 we define what we mean by
 time-consistency for processes of monetary risk measures
 and show how the time-consistency property of processes of
 monetary risk measures translates into a condition on
 processes of acceptance sets. For processes of coherent and
 convex monetary risk measures that can be represented
 with sigma-additive linear functionals we give necessary and sufficient
 conditions for time-consistency in terms of the representing
 sigma-additive linear functionals. In order to do this we define a concatenation
 operation for adapted increasing processes of integrable variation.
 The concept of m-stability for probability measures can be viewed as a special
 case of stability under concatenation. In Section 5 we discuss special
 cases and examples of time-consistent processes of monetary risk
 measures for discrete-time value processes.

 \section{The setup and notation}

 Throughout the paper $(\Omega, {\cal F}, ({\cal F}_t)_{t \in \mathbb{N}},P)$
 is a filtered probability space with ${\cal F}_0 = \crl{\emptyset, \Omega}$.
 All equalities and inequalities between random variables or
 stochastic processes are understood in the $P$-almost sure sense.
 For instance, if $(X_t)_{t \in \mathbb{N}}$ and $(Y_t)_{t \in \mathbb{N}}$
 are two stochastic processes, we mean by $X \ge Y$ that
 for $P$-almost all $\omega \in \Omega$, $X_t(\omega) \ge Y_t(\omega)$
 for all $t \in \mathbb{N}$. Also, equalities and inclusions
 between sets in ${\cal F}$ are understood in the $P$-almost sure sense.
 By ${\cal R}^0$ we denote the space of all adapted stochastic processes
 $(X_t)_{t \in \mathbb{N}}$ on $(\Omega, {\cal F}, ({\cal F}_t)_{t \in \mathbb{N}},
 P)$, where we identify two processes $X$ and $Y$ if $X = Y$.
 The two subspaces ${\cal R}^{\infty}$ and ${\cal A}^1$ of ${\cal R}^0$ are given by
 $$
 {\cal R}^{\infty} :=
 \crl{X \in {\cal R}^0 \mid \N{X}_{{\cal R}^{\infty}} < \infty} \, ,
 $$
 where
 $$
 \N{X}_{{\cal R}^{\infty}} :=
 \inf \crl{m \in \mathbb{R} \mid  \, \sup_{t \in \mathbb{N}} \abs{X_t} \le m}
 $$
 and
 $$
 {\cal A}^1 := \crl{a \in {\cal R}^0 \mid
 \N{a}_{{\cal A}^1} < \infty} \, ,
 $$
 where
 $$
 a_{-1} := 0 \, , \quad \Delta a_t := a_t - a_{t-1} \, , \mbox{ for }
 t \in \mathbb{N} \, , \quad \mbox{and} \quad
 \N{a}_{{\cal A}^1} :=
 \E{\sum_{t \in \mathbb{N}} \abs{\Delta a_t}} \, .
 $$
 The set ${\cal A}^1_+$ is given by
 $$
 {\cal A}^1_+ := \crl{a \in {\cal A}^1 \mid \Delta a_t \ge 0 \mbox{ for all }
 t \in \mathbb{N}} \, ,
 $$
 and the bilinear form $\ang{.,.}$ on ${\cal R}^{\infty}
 \times {\cal A}^1$ by
 $$
 \ang{X,a} := \E{\sum_{t \in \mathbb{N}} X_t \Delta a_t} \, .
 $$
 $\sigma({\cal R}^{\infty}, {\cal A}^1)$ denotes the coarsest
 topology on ${\cal R}^{\infty}$ such that for all $a \in {\cal A}^1$,
 $X \mapsto \ang{X,a}$ is a continuous linear functional on ${\cal
 R}^{\infty}$. $\sigma({\cal A}^1, {\cal R}^{\infty})$ denotes the coarsest
 topology on ${\cal A}^1$ such that for all $X \in {\cal R}^{\infty}$,
 $a \mapsto \ang{X,a}$ is a continuous linear functional on ${\cal
 A}^1$.

 For two $({\cal F}_t)$-stopping times $\tau$ and $\theta$ such that
 $0 \le \tau < \infty$ and $\tau \le \theta \le \infty$, we define
 the projection $\pi_{\tau, \theta} : {\cal R}^0 \to {\cal R}^0$
 by
 $$
 \pi_{\tau, \theta}(X)_t := 1_{\crl{\tau \le t}} X_{t \wedge \theta}
 \, , \quad t \in \mathbb{N} \, .
 $$
 For all $X \in {\cal R}^{\infty}$ and $a \in {\cal A}^1$, we define
 $$
 \N{X}_{\tau, \theta} := \mbox{ess\,inf }
 \crl{f \in L^{\infty}({\cal F}_{\tau}) \mid \sup_{t \in
 \mathbb{N}} \abs{\pi_{\tau, \theta}(X)_t} \le f} \, ,
 $$
 where ess\,inf denotes the essential infimum of a family of
 random variables (see for instance, Proposition VI.1.1 of Neveu,
 1972), and
 $$
 \ang{X,a}_{\tau, \theta} :=
 \E{\sum_{t \in [\tau, \theta] \cap \mathbb{N}} X_t \Delta a_t \mid {\cal
 F}_{\tau}} \, .
 $$
 The risky objects considered in this paper are elements of
 vector spaces of the form
 $$
 {\cal R}^{\infty}_{\tau, \theta} := \pi_{\tau, \theta} {\cal
 R}^{\infty} \, .
 $$
 A process $X \in {\cal R}^{\infty}_{\tau, \theta}$ is meant to
 describe the evolution of a financial value on the time interval
 $[\tau, \theta] \cap \mathbb{N}$.
 We assume that there exists a cash account where money can be
 lent to and borrowed from at the same risk-free rate and use it
 as num\'eraire, that is, all prices are expressed in
 multiples of one dollar put into the cash account at time 0.
 A monetary risk measure on ${\cal R}^{\infty}_{\tau, \theta}$ is a mapping
 $$
 \rho : {\cal R}^{\infty}_{\tau, \theta} \to L^{\infty}({\cal
 F}_{\tau}) \, ,
 $$
 assigning to a value process $X \in {\cal R}^{\infty}_{\tau, \theta}$
 a real number that can depend on the information available at the stopping
 time $\tau$ and specifies the minimal amount of money that
 has to be held in the cash account to make $X$ acceptable at time $\tau$.
 By our choice of the num\'eraire, the infusion of an amount
 of money $m$ at time $\tau$ transforms a value process $X \in {\cal
 R}^{\infty}_{\tau, \theta}$ into $X + m$ and reduces the risk
 of $X$ to $\rho(X) - m$. We find it more convenient to work with
 the negatives of monetary risk measures.
 If $\rho$ is a monetary risk measure on ${\cal R}^{\infty}_{\tau, \theta}$,
 we call $\phi = - \rho$ the monetary utility functional corresponding
 to $\rho$. $\phi(X)$ can then be viewed as a risk adjusted value of a process
 $X \in {\cal R}^{\infty}_{\tau, \theta}$ at time $\tau$.

 For the representation of conditioned convex monetary and
 coherent risk measures we will need the following subsets of
 ${\cal A}^1$:
 $$
 {\cal A}^1_{\tau, \theta} := \pi_{\tau, \theta} {\cal A}^1 \, ,
 \quad ({\cal A}^1_{\tau, \theta})_+ :=
 \pi_{\tau, \theta} {\cal A}^1_+ \quad
 \mbox{and} \quad {\cal D}_{\tau, \theta} :=
 \crl{a \in ({\cal A}^1_{\tau, \theta})_+ \mid
 \ang{1,a}_{\tau, \theta} = 1} \, .
 $$

 \section{Conditional monetary utility functionals}

 In all of Section 3, $\tau$ and $\theta$ are two fixed
 $({\cal F}_t)$-stopping times such that
 $0 \le \tau < \infty$ and $\tau \le \theta \le \infty$.

 \subsection{Basic definitions and easy properties}

 \begin{definition} \label{def1}
 We call a mapping $\phi : {\cal R}_{\tau, \theta}^{\infty} \to L^{\infty}({\cal
 F}_{\tau})$ a monetary utility functional on ${\cal R}^{\infty}_{\tau, \theta}$ if it has
 the following three properties:\\[2mm]
 {\rm (0)} $\phi(1_A X) = 1_A \phi(X)$ for all
 $X \in {\cal R}^{\infty}_{\tau, \theta}$ and $A \in {\cal
 F}_{\tau}$\\[2mm]
 {\rm (1)} $\phi(X) \le \phi(Y)$ for all $X,Y \in {\cal R}^{\infty}_{\tau, \theta}$
 such that $X \le Y$\\[2mm]
 {\rm (2)} $\phi(X + m 1_{[{\tau},\infty)}) = \phi(X) + m$ for all $X \in {\cal
 R}_{\tau}^{\infty}$ and $m \in L^{\infty}({\cal
 F}_{\tau})$\\[2mm]
 We call a monetary utility functional $\phi$ on ${\cal
 R}^{\infty}_{\tau, \theta}$ a concave monetary utility functional
 if\\
 {\rm (3)} $\phi(\lambda X + (1-\lambda) Y) \ge
 \lambda \phi(X) + (1- \lambda) \phi(Y)$ for all
 $X,Y \in {\cal R}_{\tau, \theta}^{\infty}$ and
 $\lambda \in L^{\infty}({\cal F}_{\tau})$ such that
 $0 \le \lambda \le 1$\\[2mm]
 We call a concave monetary utility functional $\phi$ on
 ${\cal R}^{\infty}_{\tau, \theta}$ a coherent utility functional if\\
 {\rm (4)} $\phi(\lambda X) = \lambda \phi(X)$
 for all $X \in {\cal R}_{\tau,\theta}^{\infty}$ and
 $\lambda \in L^{\infty}_+({\cal F}_{\tau}) :=
 \crl{f \in L^{\infty}({\cal F}_{\tau}) \mid f \ge 0}$.\\[2mm]
 For a monetary utility functional $\phi$ on ${\cal
 R}^{\infty}_{\tau, \theta}$ and $X \in {\cal R}^{\infty}$, we define
 $\phi(X) := \phi \circ \pi_{\tau, \theta}(X)$.\\[2mm]
 A monetary risk measure on ${\cal R}^{\infty}_{\tau, \theta}$
 is a mapping $\rho : {\cal R}^{\infty}_{\tau, \theta} \to
 L^{\infty}({\cal F}_{\tau})$ such that $-\rho$ is a monetary
 utility functional on ${\cal R}^{\infty}_{\tau, \theta}$.
 $\rho$ is a convex monetary risk measure if $- \rho$ is a concave
 monetary utility functional and a coherent risk measure
 if $- \rho$ is a coherent utility functional.
 \end{definition}

 \begin{Remarks} $\mbox{}$\\
 {\bf 1.} It follows from condition (0) of Definition \ref{def1} that
 $\phi(0) = 0$ for every monetary utility functional $\phi$ on
 ${\cal R}^{\infty}_{\tau, \theta}$. This normalization is convenient
 for the purposes of this paper.
 Differently normalized monetary utility functionals on ${\cal
 R}^{\infty}_{\tau, \theta}$ can be obtained by the addition of an
 ${\cal F}_{\tau}$-measurable random variable.\\
 {\bf 2.}
 It follows from (1) and (2) of Definition \ref{def1} that a monetary utility functional
 $\phi$ on ${\cal R}^{\infty}_{\tau, \theta}$ satisfies the
 following continuity condition:\\[2mm]
 (c) $\abs{\phi(X) - \phi(Y)} \le \N{X-Y}_{\tau, \theta}$,
 for all $X,Y \in {\cal R}^{\infty}_{\tau, \theta}$.\\[2mm]
 {\bf 3.} We call the property (3) of Definition \ref{def1}
 ${\cal F}_{\tau}$-concavity.\\
 {\bf 4.}
 A mapping $\phi : {\cal R}_{\tau, \theta}^{\infty} \to L^{\infty}({\cal F}_{\tau})$
 is a coherent utility functional on ${\cal R}^{\infty}_{\tau, \theta}$
 if and only if it satisfies (1), (2) and (4) of Definition \ref{def1} and\\[2mm]
 {\rm (3')} $\phi(X + Y) \ge \phi(X) + \phi(Y)$ for all
 $X,Y \in {\cal R}_{\tau, \theta}^{\infty}$.
 \end{Remarks}

 \begin{definition}
 The acceptance set ${\cal C}_{\phi}$ of a monetary utility functional
 $\phi$ on ${\cal R}^{\infty}_{\tau, \theta}$ is given by
 $$
 {\cal C}_{\phi} :=
 \crl{X \in {\cal R}^{\infty}_{\tau,\theta} \mid \phi(X) \ge 0} \,
 .
 $$
 \end{definition}
 \begin{proposition} \label{i-III}
 The acceptance set ${\cal C}_{\phi}$ of a monetary utility functional
 $\phi$ on ${\cal R}^{\infty}_{\tau, \theta}$
 satisfies the following properties:\\[2mm]
 {\rm (i)} ${\rm ess\,inf} \crl{f \in L^{\infty}({\cal F}_{\tau}) \mid
 f 1_{[\tau, \infty)} \in {\cal C}_{\phi}} = 0$.\\[2mm]
 {\rm (ii)} $1_A X + 1_{A^c} Y \in {\cal C}_{\phi}$ for all
 $X,Y \in {\cal C}_{\phi}$ and $A \in {\cal F}_{\tau}$.\\[2mm]
 {\rm (I)} $X \in {\cal C}_{\phi}$, $Y \in {\cal R}^{\infty}_{\tau, \theta}$,
 $X \le Y$ $\Rightarrow$ $Y \in {\cal C}_{\phi}$\\[2mm]
 {\rm (C)} $(X^n)_{n \in \mathbb{N}} \subset {\cal C}_{\phi}$,
 $X \in {\cal R}^{\infty}_{\tau, \theta}$,
 $\N{X^n - X}_{\tau, \theta} \stackrel{\mbox{\rm a.s.}}{\to} 0$ $\Rightarrow$ $X \in {\cal
 C}_{\phi}$.\\[2mm]
 If $\phi$ is a concave monetary utility functional, then\\[2mm]
 {\rm (II)} $\lambda X + (1-\lambda) Y \in {\cal C}_{\phi}$ for all
 $X,Y \in {\cal C}_{\phi}$ and $\lambda \in L^{\infty}({\cal F}_{\tau})$
 such that $0 \le \lambda \le 1$.\\[2mm]
 If $\phi$ is a coherent utility functional, then\\[2mm]
 {\rm (II')} $X + Y \in {\cal C}_{\phi}$ for all
 $X,Y \in {\cal C}_{\phi}$ and\\[2mm]
 {\rm (III)} $\lambda X$ for all $X \in {\cal C}_{\phi}$ and
 $\lambda \in L^{\infty}_+({\cal F}_{\tau})$.
 \end{proposition}
 \begin{proof}
 (i): It follows from the definition of ${\cal C}_{\phi}$ and
 (0) and (2) of Definition \ref{def1} that
 \beas
 && {\rm ess\,inf} \crl{f \in L^{\infty}({\cal F}_{\tau})
 \mid f 1_{[\tau, \infty)} \in {\cal C}_{\phi}}
 = {\rm ess\,inf} \crl{f \in L^{\infty}({\cal F}_{\tau})
 \mid \phi(f 1_{[\tau, \infty)}) \ge
 0}\\
 &=& {\rm ess\,inf} \crl{f \in L^{\infty}({\cal F}_{\tau})
 \mid \phi(0) + f \ge 0}
 = {\rm ess\,inf} \crl{f \in L^{\infty}({\cal F}_{\tau})
 \mid f \ge 0} = 0 \, .
 \eeas
 (ii) follows directly from (0) of Definition \ref{def1}.\\
 (I) follows from (1) of Definition \ref{def1}.\\
 (C): Let $(X^n)_{n \in \mathbb{N}}$ be a sequence in ${\cal C}_{\phi}$ and
 $X \in {\cal R}^{\infty}_{\tau, \theta}$ such that
 $\N{X^n - X}_{\tau, \theta} \stackrel{\mbox{a.s.}}{\to} 0$. It follows from (c) that
 $$
 \phi(X) \ge \phi(X^n) - \N{X^n - X}_{\tau, \theta} \, ,
 $$
 for all $n \in \mathbb{N}$. Hence, $\phi(X) \ge 0$. The remaining statements of
 the proposition are obvious.
 \end{proof}

 \begin{Remarks} $\mbox{}$\\
 {\bf 1.}
 We call a subset of ${\cal R}^{\infty}_{\tau, \theta}$ that
 satisfies condition (II) of Proposition \ref{i-III}
 ${\cal F}_{\tau}$-convex.
 {\bf 2.}
 Let ${\cal C}_{\phi}$ be the acceptance set of a monetary utility
 functional $\phi$ on ${\cal R}^{\infty}_{\tau, \theta}$.
 It can be deduced from (0) of Definition \ref{def1} or,
 alternatively, from (ii) and (C) of Proposition \ref{i-III} that
 $\sum_{n \in \mathbb{N}} 1_{A_n} X^n \in {\cal C}_{\phi}$ for every sequence
 $(X^n)_{n \in \mathbb{N}}$ in ${\cal C}_{\phi}$ and each sequence
 $(A_n)_{n \in \mathbb{N}}$ of disjoint events in ${\cal F}_{\tau}$.
 \end{Remarks}

 \begin{definition}
 If ${\cal C}$ is a subset of ${\cal R}^{\infty}_{\tau,
 \theta}$, we define for all $X \in {\cal
 R}^{\infty}_{\tau,\theta}$,
 $$
 \phi_{\cal C}(X) := {\rm ess\,sup} \crl{f \in L^{\infty}({\cal F}_{\tau}) \mid
 X - f 1_{[\tau, \infty)} \in {\cal C}} \, ,
 $$
 with the convention
 $$
 {\rm ess\,sup} \, \emptyset := - \infty \, .
 $$
 \end{definition}

 \begin{Remark}
 Note that if ${\cal C}$ satisfies (ii) of Proposition
 \ref{i-III} and, for a given $X \in {\cal R}^{\infty}_{\tau, \theta}$,
 the set
 $$
 \crl{f \in L^{\infty}({\cal F}_{\tau}) \mid
 X - f 1_{[\tau, \infty)} \in {\cal C}}
 $$
 is non-empty, then it is directed upwards, and hence,
 contains an increasing sequence $(f^n)_{n \in \mathbb{N}}$ such
 that $\lim_{n \to \infty} f^n = \phi_{\cal C}(X)$ almost surely.
 \end{Remark}

 \begin{proposition}\label{schrott}
 Let $\phi$ be a monetary utility functional on ${\cal
 R}^{\infty}_{\tau, \theta}$. Then
 $\phi_{{\cal C}_{\phi}} = \phi$.
 \end{proposition}
 \begin{proof}
 For all $X \in {\cal R}^{\infty}_{\tau, \theta}$,
 \beas
 \phi_{{\cal C}_{\phi}}(X) &=& \mbox{ess\,sup} \crl{f \in L^{\infty}({\cal F}_{\tau}) \mid
 X - f 1_{[\tau,\infty)} \in {\cal C}_{\phi}}\\
 &=& \mbox{ess\,sup} \crl{f \in L^{\infty}({\cal F}_{\tau}) \mid
 \phi(X - f 1_{[\tau,\infty)}) \ge 0}\\
 &=& \mbox{ess\,sup} \crl{f \in L^{\infty}({\cal F}_{\tau}) \mid
 \phi(X) \ge  f} = \phi(X) \, .
 \eeas
 \end{proof}

 \begin{proposition}
 If ${\cal C}$ is a subset of ${\cal R}^{\infty}_{\tau, \theta}$
 that satisfies {\rm (i)}, {\rm (ii)} and {\rm (I)} of Proposition
 \ref{i-III}, then $\phi_{\cal C}$ is a monetary utility functional
 on ${\cal R}^{\infty}_{\tau, \theta}$
 and ${\cal C}_{\phi_{\cal C}}$ is the smallest subset of
 ${\cal R}^{\infty}_{\tau, \theta}$ that contains ${\cal C}$ and
 satisfies condition {\rm (C)} of Proposition \ref{i-III}.\\
 If ${\cal C}$ satisfies {\rm (i)}, {\rm (I)} and {\rm (II)} of Proposition
 \ref{i-III}, then $\phi_{\cal C}$ is a concave monetary utility
 functional on ${\cal R}^{\infty}_{\tau, \theta}$ .\\
 If ${\cal C}$ satisfies {\rm (i)}, {\rm (I)}, {\rm (II)} and
 {\rm (III)} or {\rm (i)}, {\rm (I)}, {\rm (II')} and
 {\rm (III)} of Proposition \ref{i-III}, then $\phi_{\cal C}$
 is a coherent utility functional on ${\cal R}^{\infty}_{\tau, \theta}$.
 \end{proposition}
 \begin{proof}
 Let $X,Y \in {\cal R}^{\infty}_{\tau, \theta}$ such that
 $X = Y$ on $A$ for some $A \in {\cal F}_{\tau}$.
 Assume that $X - f 1_{[\tau,\infty)} \in {\cal C}$ for some $f \in
 L^{\infty}({\cal F}_{\tau})$.
 If ${\cal C}$ satisfies (i) and (ii) of Proposition \ref{i-III},
 there exists an $m \in L^{\infty}({\cal
 F}_{\tau})$ such that $m 1_{[\tau,\infty)} \in {\cal C}$, and
 $$
 1_A(Y-f 1_{[\tau,\infty)}) + 1_{A^c} m 1_{[\tau,\infty)} =
 1_A(X-f 1_{[\tau,\infty)}) + 1_{A^c} m 1_{[\tau,\infty)} \in {\cal C}
 \, .
 $$
 If ${\cal C}$ also satisfies (I) of Proposition \ref{i-III}, then
 $$
 Y - 1_A f 1_{[\tau,\infty)} + 1_{A^c} (m + \N{Y}_{\tau, \theta})
 1_{[\tau,\infty)} \in {\cal C} \, .
 $$
 Hence, $1_A \phi_{\cal C}(Y) \ge 1_A \phi_{\cal C}(X)$, and
 by symmetry, $1_A \phi_{\cal C}(Y) = 1_A \phi_{\cal C}(X)$.
 It follows that $\phi_{\cal C}(1_A X) = 1_A \phi_{\cal C}(X)$ for all
 $X \in {\cal R}^{\infty}_{\tau,\theta}$ and $A \in {\cal F}_{\tau}$.
 (1) of Definition \ref{def1} follows from (I) of Proposition
 \ref{i-III}. (2) of Definition \ref{def1} follows directly from the
 construction of $\phi_{\cal C}$. By Proposition \ref{i-III},
 ${\cal C}_{\phi_{\cal C}}$ satisfies condition (C) of Proposition
 \ref{i-III}, and it obviously contains ${\cal C}$.
 On the other hand, if $X \in {\cal C}_{\phi_{\cal C}}$, then there exists
 an increasing sequence $(f^n)_{n \in \mathbb{N}}$ in $L^{\infty}({\cal
 F}_{\tau})$ such that $X - f^n 1_{[\tau,\infty)} \in {\cal C}$
 and $f^n \stackrel{\mbox{a.s.}}{\to} \phi_{\cal C}(X) \ge 0$. Set $g^n := f^n
 \wedge 0$. Then, by (I) of Proposition \ref{i-III},
 $X - g^n 1_{[\tau,\infty)} \in {\cal C}$, and
 $g^n \to 0$ almost surely. Hence, ${\cal C}_{\phi_{\cal C}}$
 is the smallest subset of ${\cal R}^{\infty}_{\tau, \theta}$
 that satisfies condition (C) of Proposition \ref{i-III} and contains ${\cal C}$.
 The rest of the statements are obvious.
 \end{proof}

 \subsection{Representations for conditional concave monetary and coherent utility functionals}

 \begin{definition}
 We say a concave monetary utility functional $\phi$ on ${\cal
 R}^{\infty}_{\tau, \theta}$ is continuous for bounded decreasing sequences if
 $$
 \lim_{n \to \infty} \phi(X^{n}) = \phi(X) \quad \mbox{almost surely}
 $$
 for every decreasing sequence $(X^{n})_{n \in \mathbb{N}}$ in
 ${\cal R}^{\infty}_{\tau, \theta}$ and $X \in {\cal R}^{\infty}_{\tau, \theta}$
 such that
 $$
 X^n_t \stackrel{\mbox{\rm a.s.}}{\to} X_t \quad \mbox{for all }
 t \in \mathbb{N} \, .
 $$
 \end{definition}

 \begin{lemma} \label{KS}
 Let $\phi$ be a concave monetary utility functional on ${\cal
 R}^{\infty}_{\tau, \theta}$ that is continuous for bounded
 decreasing sequences. Then the corresponding acceptance set
 ${\cal C}_{\phi}$ is $\sigma({\cal R}^{\infty}, {\cal A}^1)$-closed.
 \end{lemma}
 \begin{proof}
 Let $(X^{\lambda})_{\lambda \in \Lambda}$ be a net in ${\cal C}_{\phi}$ and $X \in
 {\cal R}^{\infty}_{\tau, \theta}$ such that $X^{\lambda} \to X$ in
 $\sigma({\cal R}^{\infty}, {\cal A}^1)$, and assume that
 \beq \label{neg}
 \phi(X) < 0 \quad \mbox{on } A
 \eeq
 for some $A \in {\cal F}_{\tau}$ with $P[A] >0 $. Then the map
 $\tilde{\phi} : {\cal R}^{\infty} \to \mathbb{R}$ given by
 $$
 \tilde{\phi}(X) = \frac{1}{P[A]} \E{1_A \, \phi \circ \pi_{\tau, \theta}(X)}
 \, , \quad X \in {\cal R}^{\infty} \, ,
 $$
 is a concave monetary utility functional on ${\cal R}^{\infty}$
 that is continuous for bounded decreasing sequences.
 Denote by ${\cal G}$ the sigma-algebra on
 $\Omega \times \mathbb{N}$ generated by all the sets
 $A \times \crl{t}$, $t \in \mathbb{N}$, $A \in {\cal F}_t$,
 and by $\nu$ the measure on $(\Omega, {\cal G})$ given by
 $$
 \nu(A \times \crl{t}) = 2^{-(t+1)} P[A] \, , \, t \in {\mathbb N} \,
 , \, A \in {\cal F}_t \, .
 $$
 Then ${\cal R}^{\infty} = L^{\infty}(\Omega \times \mathbb{N},
 {\cal G}, \nu)$ and ${\cal A}^1$ can be identified with
 $L^1(\Omega \times \mathbb{N}, {\cal G}, \nu)$.
 Hence, it can be deduced from the Krein--\v{S}mulian theorem that
 ${\cal C}_{\tilde{\phi}}$ is a $\sigma({\cal R}^{\infty}, {\cal
 A}^1)$-closed subset of ${\cal R}^{\infty}$ (see the proof of
 Theorem 3.2 in Delbaen (2002) or Remark 4.3 in Cheridito et al. (2004)). Since
 $(X^{\lambda})_{\lambda \in \Lambda} \subset {\cal C}_{\tilde{\phi}}$, it follows that
 $$
 \frac{1}{P[A]} \E{1_A \, \phi(X)} \ge 0 \, ,
 $$
 which contradicts \eqref{neg}.
 \end{proof}

 \begin{definition} \label{defsternhag}
 For a concave monetary utility functional $\phi$ on
 ${\cal R}^{\infty}_{\tau, \theta}$ and $a \in {\cal A}^1_{\tau, \theta}$, we define
 $$
 \phi^*(a) := {\rm ess\,inf}_{X \in {\cal R}^{\infty}_{\tau, \theta}}
 \crl{\ang{X,a}_{\tau, \theta} -
 \phi(X)}
 $$
 and
 $$
 \phi^{\#}(a) := {\rm ess\,inf}_{X \in {\cal C}_{\phi}}
 \ang{X,a}_{\tau, \theta} \, .
 $$
 \end{definition}

 \begin{Remarks} \label{remsternundhag} $\mbox{}$\\
 Let $\phi$ be a concave monetary utility
 functional $\phi$ on ${\cal R}^{\infty}_{\tau, \theta}$.\\
 {\bf 1.}
 Obviously, for all $a \in {\cal A}^1_{\tau, \theta}$,
  $\phi^*(a)$ and $\phi^{\#}(a)$ are measurable functions from
 $(\Omega, {\cal F}_{\tau})$ to $[- \infty, 0]$
 and
 $$
 \phi^*(a) \le \phi^{\#}(a) \quad \mbox{for all } a \in {\cal A}^1_{\tau, \theta} \, .
 $$
 Moreover,
 \beq \label{stern=hag}
 \phi^*(a) = \phi^{\#}(a) \quad \mbox{for all } a \in {\cal
 D}_{\tau, \theta}
 \eeq
 because
 $$
 \ang{X,a}_{\tau,\theta} - \phi(X) = \ang{X - \phi(X) 1_{[\tau, \infty)} ,a}_{\tau, \theta} \,
 , \quad \mbox{and} \quad X - \phi(X) 1_{[\tau, \infty)} \in {\cal C}_{\phi} \,
 ,
 $$
 for each $X \in {\cal R}^{\infty}_{\tau, \theta}$ and $a \in {\cal D}_{\tau,
 \theta}$.\\
 {\bf 2.}
 It can easily be checked that
 $$
 \phi^{\#}(\lambda a + (1-\lambda) b)
 \ge \lambda \phi^{\#}(a) + (1-\lambda) \phi^{\#}(b) \, ,
 $$
 for all $a,b \in {\cal A}^1_{\tau, \theta}$ and
 $\lambda \in L^{\infty}({\cal F}_{\tau})$
 such that $0 \le \lambda \le 1$, and
 \beq \label{phihaglocal}
 \phi^{\#}(\lambda a) = \lambda \phi^{\#}(a)
 \quad \mbox{for all } a \in {\cal A}^1_{\tau, \theta}
 \mbox{ and } \lambda \in L_+^{\infty}({\cal F}_{\tau}) \, .
 \eeq
 Note that it follows from \eqref{phihaglocal} that
 $$
 \phi^{\#}(1_A a + 1_{A^c} b) =
 1_A \phi^{\#}(a) + 1_{A^c} \phi^{\#}(b)
 $$
 for all $a,b \in {\cal A}^1_{\tau, \theta}$ and
 $A \in {\cal F}_{\tau}$.\\
 {\bf 3.}
 For every measurable function $m: (\Omega, {\cal F}_{\tau})
 \to [- \infty, 0]$, the set
 $$
 \crl{a \in {\cal A}^1_{\tau, \theta} \mid \phi^{\#}(a) \ge m}
 $$
 is $\sigma({\cal A}^1, {\cal R}^{\infty})$-closed.
 Indeed, let $(a^{\mu})_{\mu \in M}$ be a net in
 $\crl{a \in {\cal A}^1_{\tau, \theta} \mid \phi^{\#}(a) \ge m}$
 and $a \in {\cal A}^1$ such that $a^{\mu} \to a$ in
 $\sigma({\cal A}^1, {\cal R}^{\infty})$. Then, for all $X \in {\cal
 C}_{\phi}$, $\mu \in M$
 and $A \in {\cal F}_{\tau}$ such that $A \subset \crl{m > -\infty}$,
 $$
 \ang{1_A X, a^{\mu}}
 = \E{1_A \ang{X, a^{\mu}}_{\tau, \theta}} \ge \E{1_A m} \, .
 $$
 Hence,
 $$
 \E{1_A \ang{X,a}_{\tau, \theta}} = \ang{1_A X, a} \ge \E{1_A m} \, ,
 $$
 which shows that
 $$
 \ang{X,a}_{\tau, \theta} \ge m \, , \quad \mbox{for all }
 X \in {\cal C}_{\phi} \, ,
 $$
 and therefore $\phi^{\#}(a) \ge m$.
 \end{Remarks}

 \begin{definition} \label{gamma}
 A penalty function $\gamma$ on ${\cal D}_{\tau, \theta}$ is a
 mapping from ${\cal D}_{\tau, \theta}$ to the space
 of measurable functions $f : (\Omega, {\cal F}_{\tau})
 \to [- \infty, 0]$ with the following property:
 $$
 {\rm ess\,sup}_{a \in {\cal D}_{\tau, \theta}} \gamma(a) = 0 \, .
 $$
 We call a penalty function $\gamma$ on ${\cal D}_{\tau, \theta}$
 special if
 $$
 \gamma(1_A a + 1_{A^c} b) = 1_A \gamma(a) + 1_{A^c} \gamma(b) \,
 ,
 $$
 for all $a,b \in {\cal D}_{\tau,\theta}$ and $A \in {\cal
 F}_{\tau}$.
 \end{definition}

 \begin{theorem} \label{thmrep}
 The following are equivalent:
 \begin{enumerate}
 \item[{\rm (1)}]
 $\phi$ is a mapping defined on ${\cal R}^{\infty}_{\tau, \theta}$ that can be
 represented as
 \begin{equation} \label{repr}
 \phi(X) = {\rm ess\,inf}_{a \in {\cal D}_{\tau, \theta}}
 \crl{\ang{X,a}_{\tau, \theta} -\gamma(a)} \, ,
 \quad X \in {\cal R}^{\infty}_{\tau, \theta} \, ,
 \end{equation}
 for a penalty function $\gamma$ on ${\cal D}_{\tau, \theta}$.
 \item[{\rm (2)}] $\phi$ is a concave monetary utility functional on
 ${\cal R}^{\infty}_{\tau, \theta}$
 whose acceptance set ${\cal C}_{\phi}$ is
 $\sigma({\cal R}^{\infty},{\cal A}^{1})$-closed.
 \item[{\rm (3)}] $\phi$ is a concave monetary utility functional on
 ${\cal R}^{\infty}_{\tau, \theta}$
 that is continuous for bounded decreasing sequences.
 \end{enumerate}
 Moreover, if {\rm (1)--(3)} are satisfied, then $\phi^{\#}$ is a
 special penalty function on
 ${\cal D}_{\tau, \theta}$, $\phi^{\#}(a) \ge \gamma(a)$ for all
 $a \in {\cal D}_{\tau, \theta}$, and
 the representation {\rm (\ref{repr})} also holds with $\phi^{\#}$ instead of
 $\gamma$.
 \end{theorem}

 \begin{proof} $\mbox{}$\\
 (1) $\Rightarrow$ (3):
 If $\phi$ has a representation of the form
 \eqref{repr}, then it is obviously a concave monetary utility
 functional on ${\cal R}^{\infty}_{\tau, \theta}$.
 To show that it is continuous for bounded decreasing
 sequences, let $(X^n)_{n \in \mathbb{N}}$ be a decreasing sequence in
 ${\cal R}^{\infty}_{\tau, \theta}$ and $X \in {\cal R}^{\infty}_{\tau, \theta}$
 such that
 $$
 \lim_{n \to \infty}
 X^n_t = X_t \quad \mbox{almost surely, for all } t \in \mathbb{N} \, .
 $$
 Note that this implies that
 $$
 \lim_{n \to \infty} \ang{X^n,a}_{\tau,\theta} =
 \ang{X,a}_{\tau,\theta} \quad \mbox{almost surely, for all } a \in {\cal D}_{\tau,\theta} \, .
 $$
 By property (1) of Definition \ref{def1}, $\phi_{\tau, \theta}(X^n)$
 is decreasing in $n$. Hence, almost surely, $\lim_{n \to \infty} \phi(X^n)$
 exists and $\lim_{n \to \infty} \phi(X^n) \ge \phi(X)$.
 On the other hand, there exists a sequence $(a^k)_{k \in \mathbb{N}}$ in ${\cal D}_{\tau,
 \theta}$ such that
 $$
 \phi(X) = \inf_{k \in \mathbb{N}} \crl{\ang{X,a^k} - \gamma(a^k)} \, .
 $$
 Since
 $$
 \ang{X^n,a^k} - \gamma(a^k) \ge \phi(X^n)
 $$
 for all $k,n \in \mathbb{N}$, we have that
 $$
 \ang{X,a^k} - \gamma(a^k) =
 \lim_{n \to \infty} \crl{\ang{X^n,a^k} - \gamma(a^k)} \ge
 \lim_{n \to \infty} \phi(X^n)
 $$
 for all $k \in \mathbb{N}$, and therefore also,
 $$
 \phi(X) \ge \lim_{n \to \infty} \phi(X^n) \, .
 $$
 (3) $\Rightarrow$ (2): follows from Lemma \ref{KS}.\\
 (2) $\Rightarrow$ (1):
 By \eqref{stern=hag} and the definition of $\phi^*$,
 $$
 \phi^{\#}(a) = \phi^*(a) \le \ang{X,a}_{\tau, \theta} - \phi(X)
 $$
 for all $X \in {\cal R}^{\infty}_{\tau, \theta}$ and $a \in {\cal D}_{\tau, \theta}$.
 Hence,
 \beq \label{ungl1}
 \phi(X) \le {\rm ess\,inf}_{a \in {\cal D}_{\tau,\theta}}
 \crl{\ang{X,a}_{\tau,\theta} - \phi^{\#}(a)}
 \quad \mbox{for all } X \in {\cal R}^{\infty}_{\tau, \theta} \, .
 \eeq
 To show the reverse inequality, let $m \in L^{\infty}({\cal F}_{\tau})$
 with
 \beq \label{as}
 m \le {\rm ess\,inf}_{a \in {\cal D}_{\tau, \theta}}
 \crl{\ang{X,a}_{\tau, \theta} - \phi^{\#}(a)}
 \, ,
 \eeq
 and assume that $Y = X -m 1_{[\tau,\infty)} \notin {\cal C}_{\phi}$. Since
 ${\cal C}_{\phi}$ is a convex, $\sigma ({\cal R}^{\infty}, {\cal
 A}^1)$-closed subset of ${\cal R}^{\infty}$,
 there exists an $a \in ({\cal A}^1_{\tau, \theta})_+$ such that
 $$
 \E{ \ang{Y,a}_{\tau,\theta}} = \ang{Y,a}_{0,\infty}
 < \inf_{Z \in {\cal C}_{\phi}} \ang{Z,a}_{0,\infty} =
 \E{ {\rm ess\,inf}_{Z \in {\cal C}_{\phi}}
 \ang{Z,a}_{\tau,\theta}} \, .
 $$
 Therefore, there exists a $B \in {\cal F}_{\tau}$ with $P[B] > 0$ such that
 \beq \label{ungl}
 \ang{Y,a}_{\tau,\theta} < {\rm ess\,inf}_{Z \in {\cal C}_{\phi}}
 \ang{Z,a}_{\tau,\theta} \quad \mbox{on} \quad B \, .
 \eeq
 Note that for $A = \crl{\ang{1,a}_{\tau,\theta} = 0}$,
 $$
 1_A \abs{\ang{Z,a}_{\tau,\theta}} \le 1_A \ang{\abs{Z},a}_{\tau,\theta} \le
 1_A \N{Z}_{\tau, \theta} \ang{1, a}_{\tau,
 \theta} = 0 \quad \mbox{for all}
 \quad Z \in {\cal R}^{\infty}_{\tau,\theta} \, .
 $$
 Hence, $B \subset \crl{\ang{1,a}_{\tau,\theta} > 0}$.
 Define the process $b \in {\cal D}_{\tau,\theta}$ as follows:
 $$
 b := 1_B \frac{a}{\ang{1,a}_{\tau,\theta}} + 1_{B^c} 1_{[\tau, \infty)} \,
 .
 $$
 It follows from \eqref{ungl} that
 $$
 \ang{X,b}_{\tau,\theta} - m = \ang{Y,b}_{\tau,\theta} <
 {\rm ess\,inf}_{Z \in {\cal C}_{\phi}}
 \ang{Z,b}_{\tau,\theta} = \phi^{\#}(b) \quad \mbox{on} \quad B \, .
 $$
 This contradicts \eqref{as}.
 Hence, $X - m 1_{[\tau, \infty)} \in {\cal C}_{\phi}$, and
 therefore, $\phi(X) \ge m$
 for all $m$ satisfying \eqref{as}, which shows that
 $$
 \phi(X) \ge {\rm ess\,inf}_{a \in {\cal D}_{\tau, \theta}}
 \crl{\ang{X,a}_{\tau, \theta} - \phi^{\#}(a)}\, .
 $$
 This together with \eqref{ungl1} proves that (2) implies (1) and also that
 $\phi^{\#}$ is a penalty function on ${\cal D}_{\tau, \theta}$.
 By Remark \ref{remsternundhag}.2, $\phi^{\#}$ is special.
 If $\phi$ is a concave monetary utility functional on ${\cal R}^{\infty}_{\tau, \theta}$
 with a representation of the form \eqref{repr}, then
 $$
 \ang{X,a}_{\tau,\theta} - \phi(X) \ge \gamma(a)
 $$
 for all $X \in {\cal R}^{\infty}_{\tau, \theta}$ and
 $a \in {\cal D}_{\tau, \theta}$, which implies that
 $\phi^{\#} = \phi^* \ge \gamma$ on ${\cal D}_{\tau, \theta}$.
 \end{proof}

 \begin{corollary} \label{cohrep}
 The following are equivalent:
 \begin{enumerate}
 \item[{\rm (1)}]
 $\phi$ is a mapping defined on ${\cal R}^{\infty}_{\tau, \theta}$ that can be
 represented as
 \begin{equation} \label{rep}
 \phi(X) = {\rm ess\,inf}_{a \in {\cal Q}}
 \ang{X,a}_{\tau, \theta} \, ,
 \quad X \in {\cal R}^{\infty}_{\tau, \theta} \, ,
 \end{equation}
 for a non-empty subset ${\cal Q}$ of ${\cal D}_{\tau, \theta}$.
 \item[{\rm (2)}] $\phi$ is a coherent utility functional on
 ${\cal R}^{\infty}_{\tau, \theta}$
 whose acceptance set ${\cal C}_{\phi}$ is
 $\sigma({\cal R}^{\infty},{\cal A}^{1})$-closed.
 \item[{\rm (3)}] $\phi$ is a coherent utility functional on
 ${\cal R}^{\infty}_{\tau, \theta}$
 that is continuous for bounded decreasing sequences.
 \end{enumerate}
 Moreover, if {\rm (1)--(3)} are satisfied, then the set
 $$
 {\cal Q}^0_{\phi} := \crl{a \in {\cal D}_{\tau, \theta} \mid \phi^{\#}(a) = 0}
 $$
 is equal to the smallest $\sigma({\cal A}^1, {\cal R}^{\infty})$-closed,
 ${\cal F}_{\tau}$-convex subset of ${\cal D}_{\tau, \theta}$
 that contains ${\cal Q}$, and the representation
 {\rm (\ref{repr})} also holds with ${\cal Q}^0_{\phi}$ instead of ${\cal Q}$.
 \end{corollary}
 \begin{proof}
 If (1) holds, then it follows from Theorem \ref{thmrep} that
 $\phi$ is a concave monetary utility functional on
 ${\cal R}^{\infty}_{\tau, \theta}$ that is continuous for bounded
 decreasing sequences, and it is clear that $\phi$ is coherent.
 This shows that (1) implies (3). The implication
 (3) $\Rightarrow$ (2) follows directly from Theorem \ref{thmrep}.
 If (2) holds, then Theorem \ref{thmrep} implies that
 $\phi^{\#}$ is a special penalty function on ${\cal D}_{\tau,
 \theta}$, and
 $$
 \phi(X) = \inf_{a \in {\cal D}_{\tau, \theta}}
 \crl{\ang{X,a}_{\tau, \theta} - \phi^{\#}(a)}
 \quad \mbox{for all } X \in {\cal R}^{\infty}_{\tau, \theta} \, .
 $$
 Since $\phi^{\#}$ is special, the set
 $\crl{\phi^{\#}(a) \mid a \in {\cal D}_{\tau, \theta}}$ is directed upwards.
 Therefore, there exists a sequence $(a^k)_{k \in
 \mathbb{N}}$ in ${\cal D}_{\tau, \theta}$ such that almost
 surely,
 $$
 \phi^{\#}(a^k) \nearrow \mbox{ess\,sup}_{a \in {\cal D}_{\tau, \theta}} \phi^{\#}(a) =
 0 \, , \quad \mbox{as } k \to \infty \, .
 $$
 It can easily be deduced from the fact that $\phi$ is coherent,
 that
 $$
 \crl{\phi^{\#}(a) = 0 } \cup \crl{\phi^{\#}(a) = - \infty} = \Omega
 \quad \mbox{for all } a \in {\cal D}_{\tau, \theta} \, .
 $$
 Hence, the sets $A_k := \crl{\phi^{\#}(a^k) = 0}$ are increasing in $k$, and
 $\bigcup_{k \in \mathbb{N}} A_k = \Omega$. Therefore,
 $$
 a^* := 1_{A_0} a^0 + \sum_{k \ge 1} 1_{A_k \setminus A_{k-1}} a^k
 \in {\cal D}_{\tau, \theta} \, ,
 $$
 and it follows from Remark \ref{remsternundhag}.2 that $\phi^{\#}(a^*) =
 0$. Note that for all $a \in {\cal D}_{\tau, \theta}$,
 $$
 1_{\crl{\phi^{\#}(a) = 0}} a +
 1_{\crl{\phi^{\#}(a) = - \infty}} a^* \in {\cal Q}^0_{\phi} \, .
 $$
 This shows that
 \beq \label{repq0}
 \phi(X) = \mbox{ess\,inf}_{a \in {\cal Q}^0_{\phi}} \ang{X,a}_{\tau, \theta} \, ,
 \quad \mbox{for all } X \in {\cal R}^{\infty}_{\tau, \theta} \, .
 \eeq
 It remains to show that ${\cal Q}^0_{\phi}$ is equal to the
 $\sigma({\cal A}^1, {\cal R}^{\infty})$-closed,
 ${\cal F}_{\tau}$-convex hull $\ang{\cal Q}_{\tau}$ of ${\cal Q}$. It
 follows from Theorem \ref{thmrep} that $\phi^{\#}$ is the largest
 among all penalty functions on ${\cal D}_{\tau, \theta}$ that
 induce $\phi$. This implies ${\cal Q} \subset {\cal Q}^0_{\phi}$.
 By Remarks \ref{remsternundhag}.2
 and \ref{remsternundhag}.3, ${\cal Q}^0_{\phi}$ is
 ${\cal F}_{\tau}$-convex and $\sigma({\cal A}^1, {\cal R}^{\infty})$-closed.
 Hence, $\ang{\cal Q}_{\tau} \subset {\cal Q}^0_{\phi}$. Now, assume that
 there exists a $b \in {\cal Q}^0_{\phi} \setminus \ang{\cal Q}_{\tau}$.
 Then, it follows from the separating hyperplane theorem
 that there exists an $X \in {\cal R}^{\infty}_{\tau, \theta}$
 such that
 \beq \label{gap}
 \ang{X,b} < \inf_{a \in \ang{\cal Q}_{\tau}}
 \ang{X,a} = \E{{\rm ess\,inf}_{a \in \ang{\cal Q}_{\tau}}
 \ang{X,a}_{\tau, \theta}} = \E{{\rm ess\,inf}_{a \in {\cal Q}}
 \ang{X,a}_{\tau, \theta}} = \E{\phi(X)}
 \, .
 \eeq
 But, by \eqref{repq0},
 $$
 \ang{X,b} - \E{\phi(X)}
 = \E{\ang{X,b}_{\tau, \theta} - \phi(X)} \ge 0
 $$
 for all $b \in {\cal Q}^0_{\phi}$, which contradicts \eqref{gap}.
 Hence, ${\cal Q}^0_{\phi} \setminus \ang{\cal Q}_{\tau}$ is empty, that is,
 ${\cal Q}^0_{\phi} \subset \ang{\cal Q}_{\tau}$.
 \end{proof}

 \begin{Remark}
 Detlefsen (2003) and Scandolo (2003) give
 representation results for conditional concave monetary utility
 functionals that depend on random variables.
 Since monetary utility functionals that depend
 on random variables can be seen as special cases of
 monetary utility functionals for stochastic processes, Theorem
 \ref{thmrep} generalizes the representation results in Detlefsen (2003) and
 Scandolo (2003).
 \end{Remark}

 \subsection{Relevance}

 \begin{definition}
 Let $\phi$ be a monetary utility functional on ${\cal R}^{\infty}_{\tau, \theta}$.
 We call $\phi$ $\theta$-relevant if
 $$
 A \subset \crl{\phi(- \varepsilon 1_A 1_{[t \wedge \theta, \infty)}) < 0}
 $$
 for all $\varepsilon > 0$, $t \in \mathbb{N}$ and $A \in {\cal
 F}_{t \wedge \theta}$.
 \end{definition}

 \begin{definition}
 $$
 {\cal D}^e_{\tau, \theta} :=
 \crl{a \in {\cal D}_{\tau, \theta} \mid P\edg{\sum_{j \ge t \wedge \theta} \Delta a_j
  > 0}=1 \quad
 \mbox{for all } t \in \mathbb{N}} \, .
 $$
 \end{definition}

 \begin{Remarks} $\mbox{}$\\
 {\bf 1.}
 If $\phi$ is a $\theta$-relevant monetary utility functional
 on ${\cal R}^{\infty}_{\tau, \theta}$ and
 $\xi$ is an $({\cal F}_t)$-stopping time such that
 $\tau \le \xi \le \theta$, then, obviously, the restriction of
 $\phi$ to ${\cal R}^{\infty}_{\tau, \xi}$ is $\xi$-relevant.\\
 {\bf 2.}
 Assume that $\theta$ is finite. Then it can easily be checked that
 a monetary utility functional $\phi$ on ${\cal R}^{\infty}_{\tau, \theta}$ is
 $\theta$-relevant if and only if
 $$
 A \subset \crl{\phi(- \varepsilon 1_A 1_{[\theta, \infty)}) < 0}
 $$
 for all $\varepsilon > 0$ and $A \in {\cal F}_{\theta}$.
 Also, in this case,
 $$
 {\cal D}^e_{\tau, \theta} =
 \crl{a \in {\cal D}_{\tau, \theta} \mid P\edg{\Delta a_{\theta} >
 0}=1} \, .
 $$
 \end{Remarks}

 \begin{definition}
 For a concave monetary utility functional $\phi$ on ${\cal
 R}^{\infty}_{\tau, \theta}$ and a constant $K \ge 0$, we define
 $$
 {\cal Q}^K_{\phi} :=
 \crl{a \in {\cal D}_{\tau,\theta} \mid \phi^{\#}(a) \ge - K} \, .
 $$
 \end{definition}
 By the Remarks \ref{remsternundhag}.2 and \ref{remsternundhag}.3,
 ${\cal Q}^K_{\phi}$ is ${\cal F}_{\tau}$-convex and $\sigma( {\cal A}^1,{\cal
 R}^{\infty})$-closed for every concave monetary utility functional $\phi$ on ${\cal
 R}^{\infty}_{\tau, \theta}$ and each constant $K \ge 0$.

 \begin{proposition} \label{equiv}
 Let $\phi$ be a concave monetary utility functional on
 ${\cal R}^{\infty}_{\tau, \theta}$ that is continuous for bounded
 decreasing sequences and $\theta$-relevant. Then
 $$
 {\cal Q}^K_{\phi} \cap {\cal D}^e_{\tau,\theta} \not= \emptyset \quad
 \mbox{for all} \quad K > 0 \, .
 $$
 \end{proposition}
 \begin{proof}
 Fix $K > 0$ and $t \in \mathbb{N}$.
 For $a \in {\cal D}_{\tau, \theta}$, we denote
 $$
 e_t(a) := \sum_{j \ge t \wedge \theta} \Delta a_j \, ,
 $$
 and we define
 \beq \label{alphat}
 \alpha_t := \sup_{a \in {\cal Q}^K_{\phi}} P \edg{e_t(a) > 0} \, .
 \eeq
 Let $(a^{t,n})_{n \in \mathbb{N}}$ be a sequence in ${\cal Q}^K_{\phi}$ with
 $$
 \lim_{n \to \infty} P \edg{e_t(a^{t,n}) > 0} = \alpha_t \, .
 $$
 Since ${\cal Q}^K_{\phi}$ is convex and $\sigma( {\cal A}^1,{\cal
 R}^{\infty})$-closed,
 $$
 a^t := \sum_{n \ge 1} 2^{-n} a^{t,n} \in {\cal Q}^K_{\phi} \, ,
 $$
 and, obviously,
 $$
 P \edg{e_t(a^t) > 0} = \alpha_t \, .
 $$
 In the next step we show that $\alpha_t = 1$.
 Assume to the contrary that $\alpha_t < 1$ and denote
 $A_t := \crl{e_t(a^t) = 0}$. Since $\phi$ is $\theta$-relevant,
 $$
 A_t \subset \crl{\phi(- K 1_{A_t} 1_{[t \wedge \theta,\infty)}) < 0} \, ,
 $$
 and therefore also,
 $$
 \hat{A}_t := \bigcap_{B \in {\cal F}_{\tau} , \, A_t \subset B} B
 \, \,
 \subset \, \, \crl{\phi(- K 1_{A_t} 1_{
 [t \wedge \theta, \infty)}) < 0} \, .
 $$
 By Theorem \ref{thmrep},
 $$
 \phi(- K 1_{A_t} 1_{[t \wedge \theta, \infty)})
 = {\rm ess\,inf}_{a \in {\cal D}_{\tau, \theta}}
 \crl{\ang{-K 1_{A_t} 1_{[t \wedge \theta, \infty)} , a}_{\tau,\theta} -
 \phi^{\#}(a)} \, .
 $$
 Hence, there must exist an $a \in {\cal D}_{\tau,\theta}$ with
 $P \edg{A_t \cap \crl{e_t(a) > 0}} > 0$
 and $\phi^{\#}(a) \ge -K$ on $\hat{A}_t$.
 Then,
 $$
 b^t := 1_{\hat{A}_t} \, a + 1_{\hat{A}^c} \, a^t \in {\cal
 Q}_{\phi}^K \, , \quad
 c^t := \frac{1}{2} b^t + \frac{1}{2} a^t \in {\cal
 Q}_{\phi}^K \, ,
 $$
 and $P \edg{e_t(c^t) > 0} > P \edg{e_t(a^t) > 0} = \alpha$.
 This contradicts \eqref{alphat}. Therefore, we must have $\alpha_t = 1$ for all
 $t \in \mathbb{N}$. Finally, set
 $$
 a^* = \sum_{t \ge 1} 2^{-t} a^t \, ,
 $$
 and note that $a^* \in {\cal Q}^K_{\phi} \cap {\cal D}^e_{\tau,
 \theta}$.
 \end{proof}

 \begin{corollary} \label{correpde}
 Let $\phi$ be a concave monetary utility functional
 on ${\cal R}^{\infty}_{\tau, \theta}$ that is continuous for
 bounded decreasing sequences and $\theta$-relevant. Then
 $$
 \phi(X) = {\rm ess\,inf}_{a \in {\cal D}^e_{\tau, \theta}}
 \crl{ \ang{X,a}_{\tau, \theta} - \phi^{\#}(a)} \, ,
 \quad \mbox{for all } X \in {\cal R}^{\infty}_{\tau, \theta} \, .
 $$
 \end{corollary}
 \begin{proof}
 By Theorem \ref{thmrep},
 $$
 \phi(X) = {\rm ess\,inf}_{a \in {\cal D}_{\tau, \theta}}
 \crl{ \ang{X,a}_{\tau, \theta} - \phi^{\#}(a)} \, ,
 \quad \mbox{for all } X \in {\cal R}^{\infty}_{\tau, \theta} \, ,
 $$
 which immediately shows that
 $$
 \phi(X) \le {\rm ess\,inf}_{a \in {\cal D}^e_{\tau, \theta}}
 \crl{ \ang{X,a}_{\tau, \theta} - \phi^{\#}(a)} \, ,
 \quad \mbox{for all } X \in {\cal R}^{\infty}_{\tau, \theta} \, .
 $$
 To show the reverse inequality, we choose a $b \in {\cal D}_{\tau,\theta}$.
 It follows from Proposition \ref{equiv} that there exists a process
 $c \in {\cal Q}^1_{\phi} \cap {\cal D}^e_{\tau, \theta}$.
 Then, for all $n \ge 1$,
 $$
 b^n := (1 - \frac{1}{n})b + \frac{1}{n} c \in {\cal D}^e_{\tau,\theta}
 \, ,
 $$
 $$
 \lim_{n \to \infty} \ang{X,b^n}_{\tau,\theta} =
 \lim_{n \to \infty} \crl{(1 - \frac{1}{n}) \ang{X,b}_{\tau, \theta} + \frac{1}{n}
 \ang{X,c}_{\tau,\theta}} = \ang{X,b}_{\tau,\theta} \quad
 \mbox{almost surely,}
 $$
 and
 \beas
 \phi^{\#}(b^n) &=& {\rm ess\,inf}_{X \in {\cal C}_{\phi}}
 \ang{X,b^n}_{\tau,\theta} \ge (1 -\frac{1}{n}) {\rm ess\,inf}_{X \in {\cal C}_{\phi}}
 \ang{X,b}_{\tau,\theta} + \frac{1}{n} {\rm ess\,inf}_{X \in {\cal C}_{\phi}}
 \ang{X,c}_{\tau,\theta}\\ &=& (1 -\frac{1}{n}) \phi^{\#}(b)
 + \frac{1}{n} \phi^{\#}(c) \quad \to \quad \phi^{\#}(b) \quad \mbox{almost surely} \,
 .
 \eeas
 This shows that
 $$
 \ang{X,b}_{\tau,\theta} - \phi^{\#}(b)
 \ge {\rm ess\,inf}_{a \in {\cal D}^e_{\tau,\theta}}
 \crl{ \ang{X,a}_{\tau,\theta} - \phi^{\#}(a)} \, ,
 $$
 and therefore,
 $$
 \phi(X) \ge {\rm ess\,inf}_{a \in {\cal D}^e_{\tau,\theta}}
 \crl{ \ang{X,a}_{\tau} - \phi^{\#}(a)} \, ,
 $$
 which completes the proof.
 \end{proof}

 \begin{corollary} \label{cohrepde}
 Let $\phi$ be a coherent utility functional
 on ${\cal R}^{\infty}_{\tau, \theta}$ that is
 continuous for bounded decreasing sequences and
 $\theta$-relevant.
 Then
 $$
 \phi(X) = {\rm ess\,inf}_{a \in {\cal Q}^e_{\phi}} \;
 \ang{X,a}_{\tau, \theta} \, , \quad X \in
 {\cal R}^{\infty}_{\tau, \theta} \, ,
 $$
 where ${\cal Q}^e_{\phi} := \crl{a \in {\cal D}^e_{\tau, \theta} \mid
 \phi^{\#}(a) = 0}$.
 \end{corollary}
 \begin{proof}
 This corollary can either be deduced from Corollary \ref{cohrep}
 and Proposition \ref{equiv}
 like Corollary \ref{correpde} from Theorem \ref{thmrep} and
 Proposition \ref{equiv} or from Corollary \ref{correpde}
 with the arguments used in the proof of the implication
 (2) $\Rightarrow$ (1) of Corollary \ref{cohrep}.
 \end{proof}

 \section{Processes of monetary utility functionals and acceptance sets}

 \begin{definition}
 Let $S \in \mathbb{N}$ and $T \in {\mathbb N} \cup \{\infty\}$
 such that $S \le T$. Assume that for all $t \in [S,T] \cap \mathbb{N}$,
 $\phi_{t,T}$ is a monetary utility functional on
 ${\cal R}_{t,T}^{\infty}$ with acceptance set ${\cal C}_{t,T}$.
 Then we call $(\phi_{t,T})_{t\in [S,T] \cap \mathbb{N}}$ a monetary utility
 functional process and $({\cal C}_{t,T})_{t \in [S,T] \cap \mathbb{N}}$
 an acceptance set process. We call $(\phi_{t,T})_{t \in [S,T] \cap \mathbb{N}}$
 a relevant monetary utility functional process if
 all $\phi_{t,T}$ are $T$-relevant. We call $(\phi_{t,T})_{t \in [S,T] \cap \mathbb{N}}$ a concave
 monetary utility process if $\phi_{t,T}$ is a concave monetary
 utility functional on ${\cal R}_{t,T}^{\infty}$ for all
 $t \in [S,T] \cap \mathbb{N}$. If every
 $\phi_{t,T}$ is coherent, we call $(\phi_{t,T})_{t \in [S,T] \cap \mathbb{N}}$
 a coherent utility functional process.
 \end{definition}

 \begin{definition} \label{tautheta}
 Let $S \in \mathbb{N}$ and $T \in {\mathbb N} \cup \{\infty\}$ such that $S \le T$.
 Let $(\phi_{t,T})_{t \in [S,T] \cap \mathbb{N}}$ be a monetary utility functional
 process with corresponding acceptance set process $({\cal C}_{t,T})_{t \in [S,T] \cap
 \mathbb{N}}$. Let $\tau$ and $\theta$ be two $({\cal
 F}_t)$-stopping times such that $\tau$ is finite (i.e. $\tau < \infty$) and
 $S \le \tau \le \theta \le T$.
 Then we define the mapping $\phi_{\tau, \theta} : {\cal R}^{\infty}_{\tau, \theta}
 \to L^{\infty}({\cal F}_{\tau})$ by
 \beq \label{defphitau}
 \phi_{\tau, \theta}(X) :=
 \sum_{t \in [S,T] \cap \mathbb{N}} \phi_{t,T}(1_{\crl{\tau=t}}
 X) \, ,
 \eeq
 and the set ${\cal C}_{\tau, \theta} \subset {\cal R}^{\infty}_{\tau, \theta}$ by
 \beq \label{defctau}
 {\cal C}_{\tau, \theta} := \crl{X \in {\cal R}_{\tau,\theta}^{\infty}
 \mid 1_{\crl{\tau=t}} X
 \in {\cal C}_{t,T} \mbox{ for all $t \in [S,T] \cap \mathbb N$}} \, .
 \eeq
 \end{definition}
 It can easily be checked that $\phi_{\tau, \theta}$ defined by
 \eqref{defphitau} is a monetary utility functional on
 ${\cal R}^{\infty}_{\tau, \theta}$ and that the set
 ${\cal C}_{\tau,\theta}$ given in \eqref{defctau} is
 the acceptance set of $\phi_{\tau, \theta}$. Moreover,
 if $(\phi_{t,T})_{t \in [S,T] \cap
 \mathbb{N}}$ is a concave monetary utility functional process,
 then $\phi_{\tau,\theta}$ is a concave monetary utility functional on
 ${\cal R}^{\infty}_{\tau, \theta}$. If $(\phi_{t,T})_{t \in [S,T] \cap
 \mathbb{N}}$ is coherent, then so is $\phi_{\tau,\theta}$.

 \subsection{Time-consistency}

 \begin{definition} \label{deftimecons}
 Let $S \in \mathbb{N}$ and $T \in {\mathbb N} \cup \{\infty\}$ such that $S \le T$.
 We call a monetary utility functional process
 $(\phi_{t,T})_{t \in [S,T] \cap \mathbb{N}}$ time-consistent if
 $$
 \phi_{t,T}(X) =
 \phi_{t, T}(X 1_{[t, \theta)} + \phi_{\theta,T}(X) 1_{[\theta,
 \infty)})
 $$
 for each $t \in [S,T] \cap \mathbb{N}$, every finite $({\cal
 F}_t)$-stopping time $\theta$ such that $t \le \theta \le T$ and
 all processes $X \in {\cal R}^{\infty}_{t,T}$.
 \end{definition}
 \begin{Remarks} \label{remtimeconstau}
 $\mbox{}$\\ {\bf 1.} Let $S \in \mathbb{N}$ and $T \in {\mathbb N} \cup \{\infty\}$
 such that $S \le T$. It is easy to see that a monetary utility
 functional process $(\phi_{t,T})_{t \in [S,T] \cap \mathbb{N}}$
 is time-consistent, if and only if
 $$
 \phi_{t, T}(X) \le \phi_{t, T}(Y) \, ,
 $$
 for each $t \in [S,T] \cap \mathbb{N}$ and all processes
 $X, Y \in {\cal R}^{\infty}_{t,T}$ such that
 $$
 X 1_{[t, \theta)} \le Y 1_{[t, \theta)} \quad
 \mbox{and} \quad
 \phi_{\theta, T}(X) \le \phi_{\theta,T}(Y) \, ,
 $$
 for some finite $({\cal F}_t)$-stopping time $\theta$ with $t \le
 \theta \le T$.\\
 {\bf 2.} Let $S \in \mathbb{N}$ and $T \in {\mathbb N} \cup \{\infty\}$ such that $S \le T$
 and $(\phi_{t,T})_{t \in [S,T] \cap \mathbb{N}}$ a time-consistent
 monetary utility functional process. Then it can easily be seen from
 Definition \ref{tautheta} that
 $$
 \phi_{\tau,T}(X) =
 \phi_{\tau, T}(X 1_{[\tau, \theta)} + \phi_{\theta,T}(X) 1_{[\theta, \infty)})
 $$
 for every pair of finite $({\cal F}_t)$-stopping times $\tau$ and
 $\theta$ such that $S \le \tau \le \theta \le T$ and
 all processes $X \in {\cal R}^{\infty}_{\tau,T}$.
 \end{Remarks}

 \begin{proposition} \label{proponestep}
 Let $S,T \in \mathbb{N}$ such that $S \le T$ and $(\phi_{t,T})_{t=S}^T$
 a monetary utility functional process that satisfies
 \beq \label{onestep}
 \phi_{t,T}(X) = \phi_{t,T}(X 1_{\crl{t}} + \phi_{t+1,T}(X)
 1_{[t+1, \infty)}) \,
 \eeq
 for all $t = S, \dots, T-1$ and $X \in {\cal R}^{\infty}_{t,T}$.
 Then $(\phi_{t,T})_{t=S}^T$ is time-consistent.
 \end{proposition}
 \begin{proof}
 For $t \in [S,T] \cap \mathbb{N}$, an $({\cal F}_t)$-stopping time $\theta$
 such that $t \le \theta \le T$ and a process $X \in {\cal R}^{\infty}_{t,T}$, we denote
 $Y = X 1_{[t,\theta)} + \phi_{\theta,T}(X) 1_{[\theta,\infty)}$
 and show
 \beq \label{induction}
 \phi_{t,T}(X) = \phi_{t,T}(Y)
 \eeq
 by induction. For $t=T$, \eqref{induction} is obvious. If $t \le
 T-1$, we assume that
 $$
 \phi_{t+1,T}(Z) = \phi_{t+1,T}(Z 1_{[t+1,\xi)} + \phi_{\xi,T}(Z)
 1_{[\xi, \infty)}) \, ,
 $$
 for every $({\cal F}_t)$-stopping time $\xi$ such that $t+1 \le
 \xi \le T$ and all $Z \in {\cal R}^{\infty}_{t+1,T}$. Then
 $$
 1_{\crl{\theta \ge t +1}} \phi_{t+1}(X)
 = \phi_{t+1}(1_{\crl{\theta \ge t +1}} X)
 = \phi_{t+1}(1_{\crl{\theta \ge t +1}} Y)
 = 1_{\crl{\theta \ge t +1}} \phi_{t+1}(Y) \, .
 $$
 Hence, it follows from the assumption \eqref{onestep} that
 \beas
 \phi_{t,T}(Y) &=&
 \phi_{t,T} \brak{1_{\crl{\theta = t}} \phi_{t,T}(X) 1_{[t,\infty)}
 + 1_{\crl{\theta \ge t +1}} Y}\\
 &=& 1_{\crl{\theta = t}} \phi_{t,T}(X)
 + 1_{\crl{\theta \ge t +1}} \phi_{t,T} (Y)\\
 &=& 1_{\crl{\theta = t}} \phi_{t,T}(X)
 + 1_{\crl{\theta \ge t +1}} \phi_{t,T}(Y 1_{\crl{t}} +
 \phi_{t+1}(Y) 1_{[t+1, \infty)})\\
 &=& 1_{\crl{\theta = t}} \phi_{t,T}(X)
 + 1_{\crl{\theta \ge t +1}} \phi_{t,T}(X 1_{\crl{t}} +
 \phi_{t+1}(X) 1_{[t+1, \infty)})\\
 &=& 1_{\crl{\theta = t}} \phi_{t,T}(X)
 + 1_{\crl{\theta \ge t +1}} \phi_{t,T}(X)\\
 &=& \phi_{t,T}(X) \, .
 \eeas
 \end{proof}

 \begin{proposition} \label{decomp}
 Let $S \in \mathbb{N}$, $T \in {\mathbb N} \cup \{\infty\}$ and
 $\tau, \theta$ finite $({\cal F}_t)$-stopping times
 such that $S \le \tau \le \theta \le T$.
 For a monetary utility functional process
 $(\phi_{t,T})_{t \in [S,T] \cap \mathbb{N}}$ with corresponding acceptance
 set process $({\cal C}_{t,T})_{t \in [S,T] \cap \mathbb{N}}$ the following
 two conditions are equivalent:
 \begin{itemize}
 \item[\rm (1)]
 $\phi_{\tau, T}(X) = \phi_{\tau, T}(X 1_{[\tau, \theta)} +
 \phi_{\theta,T}(X ) 1_{[\theta,\infty)})$
 for all $X \in {\cal R}^{\infty}_{\tau,T}$.
 \item[\rm (2)]
 ${\cal C}_{\tau,T} = {\cal C}_{\tau,\theta} + {\cal C}_{\theta,T}$
 \end{itemize}
 \end{proposition}
 \begin{proof}

 (1) $\Rightarrow$ (2):\\
 Assume $Y \in {\cal C}_{\tau,\theta}$ and $Z \in {\cal C}_{\theta,T}$.
 Then $X = Y+Z \in {\cal R}^{\infty}_{\tau,T}$,
 $X 1_{[\tau,\theta)} = Y 1_{[\tau,\theta)}$ and
 $\phi_{\theta,T}(X) = Y_{\theta} + \phi_{\theta,T}(Z) \ge
 Y_{\theta}$.
 Therefore,
 $$
 \phi_{\tau,T}(X) = \phi_{\tau, T}(X 1_{[\tau, \theta)} +
 \phi_{\theta,T}(X ) 1_{[\theta,\infty)}) \ge \phi_{\tau,T}(Y) \ge
 0 \, .
 $$
 This shows that ${\cal C}_{\tau,\theta} + {\cal C}_{\theta,T} \subset {\cal
 C}_{\tau,T}$. To show ${\cal C}_{\tau,T} \subset {\cal C}_{\tau,\theta} + {\cal
 C}_{\theta,T}$, let $X \in {\cal C}_{\tau,T}$ and set
 $Z := (X - \phi_{\theta,T}(X)) 1_{[\theta,\infty)}$ and $Y :=
 X - Z = X 1_{[\tau, \theta)} + \phi_{\theta,T}(X) 1_{[\theta,\infty)}$.
 It follows directly from the translation invariance of
 $\phi_{\theta,T}$ that $Z \in {\cal C}_{\theta,T}$. Moreover,
 $\phi_{\tau,T}(Y) = \phi_{\tau,T}(X) \ge 0$, which shows that
 $Y \in {\cal C}_{\tau, \theta}$.

 (2) $\Rightarrow$ (1):\\
 Let $X \in {\cal R}^{\infty}_{\tau,T}$ and
 $f \in L^{\infty}({\cal F}_{\tau})$ such that
 $X - f 1_{[\tau,\infty)} \in {\cal C}_{\tau,T}$. Since
 $$
 \phi_{\theta,T}(X ) = \mbox{ess\,sup}
 \crl{g \in L^{\infty}({\cal F}_{\theta}) \mid
 (X - g) 1_{[\theta,\infty)} \in {\cal C}_{\theta,T}}
 $$
 and
 $$
 {\cal C}_{\tau,T} \subset {\cal C}_{\tau,\theta} + {\cal
 C}_{\theta,T} \, ,
 $$
 the process
 $$
 X - f 1_{[\tau,\infty)} - (X - \phi_{\theta,T}(X)) 1_{[\theta,\infty)} =
 X 1_{[\tau, \theta)} + \phi_{\theta,T}(X) 1_{[\theta, \infty)}
 - f 1_{[\tau,\infty)}
 $$
 has to  be in ${\cal C}_{\tau,\theta}$. This shows that
 $$
 \phi_{\tau,T}(X 1_{[\tau, \theta)} + \phi_{\theta,T}(X) 1_{[\theta,
 \infty)}) \ge \phi_{\tau,T}(X) \, .
 $$
 On the other hand, if $X \in {\cal R}^{\infty}_{\tau,T}$ and
 $f \in L^{\infty}({\cal F}_{\tau})$ such that
 $$
 X 1_{[\tau, \theta)} + \phi_{\theta,T}(X) 1_{[\theta, \infty)} - f 1_{[\tau,\infty)}
 \in {\cal C}_{\tau,T} \, ,
 $$
 then also $X - f 1_{[\tau,\infty)} \in {\cal C}_{\tau,T}$ because
 $(X - \phi_{\theta,T}(X)) 1_{[\theta, \infty)} \in {\cal C}_{\theta,T}$
 and ${\cal C}_{\tau,\theta} + {\cal C}_{\theta,T} \subset {\cal
 C}_{\tau,T}$. It follows that
 $$
 \phi_{\tau,T}(X) \ge \phi_{\tau,T}(X 1_{[\tau, \theta)} + \phi_{\theta,T}(X) 1_{[\theta,
 \infty)}) \, .
 $$
 \end{proof}

 \begin{proposition} \label{cthetactau}
 Let $S \in \mathbb{N}$ and $T \in {\mathbb N} \cup \{\infty\}$ such that $S \le T$.
 Let $(\phi_{t,T})_{t \in [S,T] \cap \mathbb{N}}$ be a
 time-consistent monetary utility functional
 process with corresponding acceptance set process $({\cal C}_{t,T})_{t \in [S,T] \cap
 \mathbb{N}}$, and let $\tau$ and $\theta$ be two finite $({\cal F}_t)$-stopping times
 such that $S \le \tau \le \theta \le T$. Then\\[2mm]
 {\bf 1.}
 $1_A X \in {\cal C}_{\tau, T}$ for all $X \in
 {\cal C}_{\theta,T}$ and $A \in {\cal F}_{\theta}$.\\[2mm]
 {\bf 2.}
 If $\phi_{\tau,\theta}$ is $\theta$-relevant, and $X$ is a process in
 ${\cal R}^{\infty}_{\theta,T}$ such that
 $1_A X \in {\cal C}_{\tau,T}$ for all $A \in {\cal F}_{\theta}$,
 then $X \in {\cal C}_{\theta,T}$.\\[2mm]
 {\bf 3.}
 If $\xi$ is an $({\cal F}_t)$-stopping time such that
 $\theta \le \xi \le T$ and
 $\phi_{\tau,\xi}$ is $\xi$-relevant, then
 $\phi_{\theta,\xi}$ is $\xi$-relevant too.
 In particular, the monetary utility functional process $(\phi_{t,T})_{t \in [0,T]
 \cap \mathbb{N}}$ is relevant if and only if $\phi_{S,T}$ is $T$-relevant.
 \end{proposition}
 \begin{proof} $\mbox{}$\\
 1.
 If $X \in {\cal C}_{\theta,T}$ and $A \in {\cal F}_{\theta}$,
 then also $1_A X \in {\cal C}_{\theta,T}$.
 Obviously, $0 \in {\cal C}_{\tau,\theta}$. Hence, it follows from
 Proposition \ref{decomp} that $1_A X = 0 + 1_A X \in {\cal C}_{\tau,T}$.\\
 2.
 Assume $1_A X \in {\cal C}_{\tau,T}$ for all $A \in {\cal F}_{\theta}$
 but $X \notin {\cal C}_{\theta,T}$. Then there exists
 an $\varepsilon > 0$ such that $P[A] > 0$, where
 $A = \crl{\phi_{\theta,T}(X) \le - \varepsilon}$.
 By Proposition \ref{decomp}, there exist $Y \in {\cal C}_{\tau,\theta}$
 and $Z \in {\cal C}_{\theta,T}$ such that $1_A X = Y + Z$.
 Since $\phi_{\theta,T}(1_A X) \le - \varepsilon 1_A$,
 $Z_{\theta} \ge 1_A (X_{\theta} + \varepsilon)$ and therefore,
 $Y_{\theta} \le - \varepsilon 1_A$. But then, since $\phi_{\tau,\theta}$
 is $\theta$-relevant, $Y \notin {\cal C}_{\tau,\theta}$, which is a
 contradiction.\\
 3.
 Let $\varepsilon > 0$, $t \in \mathbb{N}$ and $A \in {\cal F}_{t
 \wedge \xi}$. Set
 $$
 B := A \cap \crl{\phi_{\theta, \xi}(-\varepsilon 1_A
 1_{[t \wedge \xi, \infty)}) = 0}
 $$
 and note that
 $$
 \phi_{\theta, \xi}(-\varepsilon 1_B 1_{[t \wedge \xi, \infty)}) = 0
 \quad \mbox{on } B \, .
 $$
 Therefore, also
 $$
 \phi_{\theta, \xi}(-\varepsilon 1_B 1_{[t \wedge \xi, \infty)}) = 0
 \quad \mbox{on } \hat{B} := \bigcap_{C \in {\cal F}_{\theta} \, ; \, B \subset C} C \, .
 $$
 Since
 $$
 1_{\hat{B}^c} \,
 \phi_{\theta, \xi}(-\varepsilon 1_B 1_{[t \wedge \xi, \infty)}) =
 \phi_{\theta, \xi}(-\varepsilon 1_{\hat{B}^c} 1_B 1_{[t \wedge \xi, \infty)}) =
 0 \, ,
 $$
 it follows that $\phi_{\theta, \xi}(-\varepsilon 1_B 1_{[t \wedge \xi, \infty)}) =
 0$. Hence, $-\varepsilon 1_B 1_{[t \wedge \xi, \infty)} \in {\cal C}_{\theta,
 \xi}$, and therefore, by statement 1, $-\varepsilon 1_B 1_{[t \wedge \xi, \infty)}
 \in {\cal C}_{\tau, \xi}$. If $\phi_{\tau, \xi}$ is
 $\xi$-relevant, then $P[B] = 0$, which shows that
 $\phi_{\theta, \xi}$ is $\xi$-relevant.
 \end{proof}

 \begin{corollary}
 Let $S \in \mathbb{N}$ and $T \in {\mathbb N} \cup \{\infty\}$ such that $S \le T$.
 Let $\phi$ be a $T$-relevant monetary utility functional on
 ${\cal R}^{\infty}_{S,T}$. Then there exists at most one
 time-consistent monetary utility process $(\phi_{t,T})_{t \in [S,T]
 \cap \mathbb{N}}$ with $\phi_{S,T} = \phi$.
 \end{corollary}
 \begin{proof}
 Let $(\phi_{t,T})_{t \in [S,T] \cap
 \mathbb{N}}$ be a time-consistent monetary utility process
 with $\phi_{S,T} = \phi$ and
 $({\cal C}_{t,T})_{t \in [S,T] \cap \mathbb{N}}$ the corresponding
 acceptance set process. By Proposition \ref{cthetactau}.3,
 $\phi_{t,T}$ is $T$-relevant for all $t \in [S,T] \cap
 \mathbb{N}$. Therefore it follows from 1. and 2. of
 Proposition \ref{cthetactau} that for all
 $t \in [S,T] \cap \mathbb{N}$, a process $X \in {\cal
 R}^{\infty}_{t,T}$ is in ${\cal C}_{t,T}$ if and only if
 $1_A X \in {\cal C}_{S,T}$ for all $A \in {\cal F}_t$.
 This shows that ${\cal C}_{t, T}$ is uniquely determined by
 the acceptance set ${\cal C}_{S,T}$ of $\phi$. Hence,
 $\phi_{t,T}$ is uniquely determined by $\phi$.
 \end{proof}

 \subsection{Consistent extension of the time horizon.}

 \begin{proposition} \label{propextension}
 Let $S \in \mathbb{N}$ and $T \in {\mathbb N} \cup \{\infty\}$ such that $S \le T$.
 Let $(\phi_{t,S})_{t \in [0,S] \cap \mathbb{N}}$ and $(\phi_{t,T})_{t \in [S,T] \cap
 \mathbb{N}}$ be two time-consistent monetary utility functional
 processes with corresponding acceptance set processes
 $({\cal C}_{t,S})_{t \in [0,S] \cap \mathbb{N}}$ and
 $({\cal C}_{t,T})_{t \in [S,T] \cap \mathbb{N}}$, respectively.
 For $t \in [0,S)$, define
 \beq \label{extphi}
 \phi_{t,T}(X) := \phi_{t,S} \brak{X 1_{[t,S)} + \phi_{S,T}(X)
 1_{[S,\infty)}} \, , \quad X \in {\cal R}^{\infty}_{t,T} \, ,
 \eeq
 and
 \beq \label{extc}
 {\cal C}_{t,T} := {\cal C}_{t,S} + {\cal C}_{S,T} \, .
 \eeq
 Then $(\phi_{t,T})_{t \in [0,T] \cap \mathbb{N}}$ is a
 time-consistent monetary utility functional process with
 corresponding acceptance set process
 $({\cal C}_{t,T})_{t \in [0,T] \cap \mathbb{N}}$.
 If $(\phi_{t,S})_{t \in [0,S] \cap \mathbb{N}}$ and $(\phi_{t,T})_{t \in [S,T] \cap
 \mathbb{N}}$ are concave monetary utility functional processes,
 then so is $(\phi_{t,T})_{t \in [0,T] \cap \mathbb{N}}$.
 If $(\phi_{t,S})_{t \in [0,S]} \cap \mathbb{N}$ and $(\phi_{t,T})_{t \in [S,T] \cap
 \mathbb{N}}$ are coherent, then $(\phi_{t,T})_{t \in [0,T] \cap \mathbb{N}}$
 is coherent too.
 \end{proposition}
 \begin{proof}
 It can easily be checked that for all $t \in [0,S)$, the
 mapping $\phi_{t,T}$ defined in \eqref{extphi} is a monetary
 utility functional on ${\cal R}^{\infty}_{t,T}$ with acceptance
 set ${\cal C}_{t,T}$ given by \eqref{extc}. Also, it is obvious
 that $\phi_{t,T}$ is a concave monetary utility functional on
 ${\cal R}^{\infty}_{t,T}$ if $\phi_{t,S}$ and $\phi_{S,T}$ are concave
 monetary utility functionals, and $\phi_{t,T}$ is
 coherent if $\phi_{t,S}$ and $\phi_{S,T}$ are coherent.
 To prove that $(\phi_{t,T})_{t \in [0,T] \cap \mathbb{N}}$
 is time-consistent, it is by Proposition \ref{decomp} enough to show that
 $$
 {\cal C}_{t,T} = {\cal C}_{t, \theta} + {\cal C}_{\theta,T} \, ,
 $$
 for all $t \in [0,S) \cap \mathbb{N}$ and every finite
 $({\cal F}_t)$-stopping time $\theta$ such that $t \le \theta \le T$.

 We first show ${\cal C}_{t, \theta} + {\cal C}_{\theta,T} \subset {\cal
 C}_{t,T}$. Let $Y \in {\cal C}_{t,\theta}$ and $Z \in {\cal C}_{\theta,T}$.
 By definition of ${\cal C}_{t,T}$, $Y$ can be decomposed into $Y =
 Y' + Y''$, where $Y' \in {\cal C}_{t,S}$ and $Y'' \in
 {\cal C}_{S,T}$. It is easy to see that $Y''$ can be chosen such
 that $Y'' = 0$ on $\crl{\theta \le S}$. Then
 $$
 Y' \in {\cal C}_{t, \theta \wedge S} \quad \mbox{and} \quad
 Y'' \in {\cal C}_{S, \theta \vee S} \,.
 $$
 Similarly, $Z = Z' + Z''$, where $Z' \in {\cal C}_{t,S}$ and
 $Z'' \in {\cal C}_{S,T}$ can be chosen such that
 $$
 Z' \in {\cal C}_{\theta \wedge S, S} \quad
 \mbox{and} \quad
 Z'' \in {\cal C}_{\theta \vee S, T} \, .
 $$
 Hence,
 $$
 Y' + Z' \in {\cal C}_{t,S} \quad \mbox{and} \quad
 Y'' + Z'' \in {\cal C}_{S,T} \, ,
 $$
 and therefore,
 $$
 Y+Z = Y' + Z' + Y'' + Z'' \in {\cal C}_{t,T} \, .
 $$

 To show ${\cal C}_{t,T} \subset {\cal C}_{t, \theta} + {\cal
 C}_{\theta,T}$, we let $X \in {\cal C}_{t,T}$. By definition of
 ${\cal C}_{t,T}$, $X = X' + X''$, where
 $X' \in {\cal C}_{t,S}$ and $X'' \in {\cal C}_{S,T}$.
 Since $(\phi_{t,S})_{t \in [0,S] \cap \mathbb{N}}$
 and $(\phi_{t,T})_{t \in [S,T] \cap \mathbb{N}}$
 are time-consistent, we get from Proposition \ref{decomp} that
 $$
 X' = Y' + Z' \quad \mbox{and} \quad X'' = Y'' + Z'' \, ,
 $$
 where $Y' \in {\cal C}_{t, \theta \wedge S}$, $Z' \in {\cal
 C}_{\theta \wedge S, S}$, $Y'' \in {\cal C}_{S, \theta \vee S}$
 and $Z'' \in {\cal C}_{\theta \vee S, T}$. Note that
 $1_{\crl{\theta > S}} Z' \in {\cal C}_{S,S}$. Hence,
 $$
 1_{\crl{\theta > S}} Z' = f 1_{\crl{\theta >S}} 1_{[S, \infty)}
 \quad \mbox{for some } f \in L^{\infty}_+({\cal F}_S) \, .
 $$
 It follows that
 $$
 1_{\crl{\theta > S}} Y'' + 1_{\crl{\theta > S}} Z' \in
 {\cal C}_{S, \theta \vee S} \cap {\cal R}^{\infty}_{t, \theta} \,
 ,
 $$
 and therefore,
 $$
 Y' + 1_{\crl{\theta > S}} Y'' + 1_{\crl{\theta > S}} Z' \in {\cal
 C}_{t, \theta} \, .
 $$
 On the other hand,
 $$
 1_{\crl{\theta \le S}} Y'' \in {\cal C}_{S,S} \, ,
 $$
 and therefore,
 $$
 1_{\crl{\theta \le S}} Y''
 = g 1_{\crl{\theta \le S}} 1_{[S, \infty)} \quad
 \mbox{for some } g \in L^{\infty}_+({\cal F}_S) \, .
 $$
 Hence,
 $$
 1_{\crl{\theta \le S}} Y'' + 1_{\crl{\theta \le S}} Z'
 \in {\cal C}_{\theta \wedge S, S} \cap {\cal R}^{\infty}_{\theta,
 T} \, ,
 $$
 and
 $$
 1_{\crl{\theta \le S}} Y'' + 1_{\crl{\theta \le S}} Z' + Z''
 \in {\cal C}_{\theta, T} \, .
 $$
 \end{proof}

 \begin{Remark}
 Let $T \in \mathbb{N}$. Note that for all
 $t = 0, \dots, T$, there exists only one
 monetary utility functional $\phi_{t,t}$
 on ${\cal R}^{\infty}_{t,t}$. It is given by
 $$
 \phi_{t,t}(m 1_{[t,\infty)}) = m \, , \quad
 \mbox{for } m \in L^{\infty}({\cal F}_t) \, ,
 $$
 and its acceptance set is
 $$
 C_{t,t} = \crl{m 1_{[t, \infty)} \mid m \in
 L^{\infty}_+(\Omega, {\cal F}_t ,P)} \, .
 $$
 Now, for every $t = 0, \dots T-1$, let
 $\phi_{t,t+1}$ be an arbitrary monetary utility functional
 on ${\cal R}^{\infty}_{t,t+1}$ with acceptance set
 ${\cal C}_{t, t+1}$. It can easily be checked that
 for all $t = 0, \dots T-1$, the monetary utility functional process
 $(\phi_{s,t+1})_{s=t}^{t+1}$ is time-consistent.
 Therefore, it follows from Proposition \ref{propextension} that an
 acceptance set process $({\cal C}_{t,T})_{t \in [0,T] \cap
 \mathbb{N}}$ corresponding to a time-consistent monetary utility
 functional process $(\phi_{t,T})_{t \in [0,T] \cap \mathbb{N}}$
 can be obtained by defining
 $$
 {\cal C}_{t,T} := {\cal C}_{t,t+1} + {\cal C}_{t+1,t+2} + \dots
 + {\cal C}_{T-1,T} \, , \quad \mbox{for all } t = 0,1, \dots,T-1 \, .
 $$
 \end{Remark}

 \subsection{Concatenation of elements in ${\cal A}^1_+$}

 \begin{definition}
 Let $a,b \in {\cal A}^1_+$, $\theta$ a finite
 $({\cal F}_t)$-stopping time and $A \in {\cal F}_{\theta}$. Then the
 concatenation $a \oplus^{\theta}_A b$ is defined by
 $$
 (a \oplus^{\theta}_A b)_t :=
 \left\{
 \begin{array}{ll}
 a_t & \mbox{on }
 \crl{t < \theta} \cup A^c \cup \crl{\ang{1,b}_{\theta, \infty} = 0}\\
 a_{\theta-1} + \frac{\ang{1,a}_{\theta,\infty}}{\ang{1,b}_{\theta,\infty}}
 \brak{b_t - b_{\theta - 1}} & \mbox{on }
 \crl{t \ge \theta} \cap A \cap \crl{\ang{1,b}_{\theta, \infty} > 0} \, ,
 \end{array}
 \right.
 $$
 where we set $a_{-1} = b_{-1} = 0$.\\[2mm]
 We say a subset ${\cal Q}$ of ${\cal A}^1_+$ is stable under
 concatenation if $a \oplus^{\theta}_A b \in {\cal Q}$ for all
 $a,b \in {\cal Q}$, every finite $({\cal F}_t)$-stopping time
 $\theta $ and all $A \in {\cal F}_{\theta}$.
 \end{definition}

 \begin{Remarks} \label{remprop} $\mbox{}$\\
 {\bf 1.} Let ${\cal Q}$ be a subset of ${\cal A}^1_+$ such that
 \beq \label{concats}
 a \oplus^s_A b \in {\cal Q} \quad \mbox{for all }
 a,b \in {\cal Q} \, , \, s \in \mathbb{N} \mbox{ and }
 A \in {\cal F}_s \, .
 \eeq
 Then
 $$
 a \oplus^{\theta}_A b \in {\cal Q} \quad \mbox{for all }
 a,b \in {\cal Q}, \mbox{ each bounded } ({\cal F}_t) \mbox{-stopping
 time } \theta, \mbox{ and } A \in \cal F_{\theta} \, ,
 $$
 and
 $$
 a \oplus^{\theta}_A b \quad
 \mbox{is in the } \N{.}_{{\cal A}^1} \mbox{-closure of }
 {\cal Q}
 $$
 for all $a,b \in {\cal Q}$, each finite $({\cal F}_t)$-stopping
 time $\theta$ and $A \in {\cal F}_{\theta}$.

 Indeed, if ${\cal Q}$ has the property \eqref{concats}, set for each
 $({\cal F}_t)$-stopping time $\theta$ and $A \in {\cal
 F}_{\theta}$, $A_n := A \cap \crl{\theta = n}$, $n \in \mathbb{N}$.
 Then all the following processes are in ${\cal Q}$:
 $$
 a^0 := a \oplus^0_{A_0} b \, , \quad a^n := a^{n-1} \oplus^n_{A_n} b \, ,
 \quad n \ge 1 \, .
 $$
 If $\theta$ is bounded, then $a^n = a \oplus^{\theta}_A b$ for
 all $n$ such that $n \ge \theta$. If $\theta$ is finite, then
 $a^n \to a \oplus^{\theta}_A b$ in $\N{.}_{\cal A}^1$, as $n \to \infty$.\\
 {\bf 2.}
 Let $\theta$ be a finite $({\cal F}_t)$-stopping time and
 $A \in {\cal F}_{\theta}$. It can easily be checked that the concatenation
 $\oplus^{\theta}_A$ has the following properties:\\
 (i) Let $a^1, a^2, b \in {\cal A}^1_+$ and $\lambda \in (0,1)$.
 Then
 $$
 (\lambda a^1 + (1 - \lambda) a^2) \oplus^{\theta}_A b
 = \lambda (a^1 \oplus^{\theta}_A b) + (1 - \lambda) (a^2 \oplus^{\theta}_A
 b)\, .
 $$
 (ii)
 Let $a,b \in {\cal A}^1_+$ and $(a^{\mu})_{\mu \in M}$
 a net in ${\cal A}^1_+$ with
 $$
 a^{\mu} \to a \quad \mbox{in }
 \sigma({\cal A}^1, {\cal R}^{\infty}) \, .
 $$
 Then
 $$
 a^{\mu} \oplus^{\theta}_A b \to a \oplus^{\theta}_A b \quad \mbox{in }
 \sigma({\cal A}^1, {\cal R}^{\infty}) \, .
 $$
 (iii)
 Let $a,b \in {\cal A}^1_+$, $B = \crl{\ang{1,b}_{\theta} >0}$
 and $(b^{\mu})_{\mu \in M}$ a net in ${\cal A}^1_+$ with
 $$
 b^{\mu} \to b \quad \mbox{in }
 \sigma({\cal A}^1, {\cal R}^{\infty}) \, .
 $$
 Then
 $$
 a \oplus^{\theta}_{A \cap B} b^{\mu} \to a \oplus^{\theta}_A b \quad \mbox{in }
 \sigma({\cal A}^1, {\cal R}^{\infty}) \, .
 $$
 \end{Remarks}

 \begin{proposition} \label{estable}
 Let ${\cal Q}$ be a non-empty subset of ${\cal A}^1_+$ and denote by
 $\ang{\cal Q}$ the smallest $\sigma({\cal A}^1, {\cal
 R}^{\infty})$-closed, convex subset of ${\cal A}^1_+$ that contains ${\cal Q}$.
 Assume that
 $$
 a \oplus^s_A b \in \ang{\cal Q} \, ,
 $$
 for all $a,b \in {\cal Q}$, $s \in \mathbb{N}$ and $A \in {\cal F}_s$.
 Then $\ang{\cal Q}$ is stable under concatenation.
 \end{proposition}
 \begin{proof}
 It follows from the properties (i) and (ii) of Remark \ref{remprop}.2 that
 $$
 a \oplus^s_A b \in \ang{\cal Q}
 $$
 for all $a \in \ang{\cal Q}$, $b \in {\cal Q}$, $s \in \mathbb{N}$ and $A \in {\cal
 F}_s$. Then, it can be shown as in Remark \ref{remprop}.1 that
 $$
 a \oplus^{\theta}_A b \in \ang{\cal Q}
 $$
 for all $a \in \ang{\cal Q}$, $b \in {\cal Q}$, every finite
 $({\cal F}_t)$-stopping time $\theta$ and $A \in {\cal F}_{\theta}$.
 Now, let $a \in \ang{\cal Q}, b^1, b^2 \in {\cal Q}$ and
 $\lambda \in (0,1)$. Set
 $$
 B_1 := \crl{\ang{1,b^1}_{\theta} > 0} \quad \mbox{and} \quad
 B_2 := \crl{\ang{1,b^2}_{\theta} > 0} \, .
 $$
 Let $C_1, \dots, C_N$ be finitely many disjoint sets in ${\cal
 F}_{\theta}$ such that $\bigcup_{n=1}^N C_n = B_1 \cap B_2$ and
 $\lambda_1, \dots, \lambda_N$ numbers in $[0,1]$.
 Then the following processes are all in $\ang{\cal Q}$:
 $$
 c^1 = \lambda_1 (a \oplus^{\theta}_{C_1} b^1) + (1 - \lambda_1)
 (a \oplus^{\theta}_{C_1} b^2) \, ,
 $$
 $$
 c^2 = \lambda_2 (c^1 \oplus^{\theta}_{C_2} b^1) + (1 - \lambda_2)
 (c^1 \oplus^{\theta}_{C_2} b^2) \, , \dots
 $$
 $$
 \dots \, , \, c^N = \lambda_N (c^{N-1} \oplus^{\theta}_{C_N} b^1) + (1 - \lambda_N)
 (c^{N-1} \oplus^{\theta}_{C_N} b^2) \, .
 $$
 Note that
 $$
 c^N_t =
 \left\{
 \begin{array}{ll}
 a_t & \mbox{ on } \crl{t < \theta} \cup B^c_1 \cup B^c_2\\
 a_{\theta - 1} +
 \lambda_n \frac{\ang{1,a}_{\theta}}{\ang{1,b^1}_{\theta}}
 (b^1_t - b^1_{\theta-1})
 + (1-\lambda_n) \frac{\ang{1,a}_{\theta}}{\ang{1,b^2}_{\theta}}
 (b^2_t - b^2_{\theta-1}) & \mbox{ on }
 \crl{t \ge \theta} \cap C_n
 \end{array}
 \right. \, .
 $$
 Hence, since $\ang{\cal Q}$ is $\sigma({\cal A}^1, {\cal R}^{\infty})$-closed,
 it also contains the process $c$ given by
 $$
 c_t =
 \left\{
 \begin{array}{ll}
 a_t & \mbox{ on } \crl{t < \theta} \cup B^c_1 \cup B^c_2\\
 a_{\theta-1} + \frac{\ang{1,a}_{\theta}}{\ang{1, \lambda b^1 + (1 - \lambda)
 b^2}_{\theta}} \edg{\lambda (b^1_t - b^1_{\theta-1}) +
 (1 - \lambda) (b^2_t - b^2_{\theta}) } & \mbox{ on }
 \crl{t \ge \theta} \cap B^1 \cap B^2
 \end{array}
 \right. \, .
 $$
 Next, notice that the processes
 $$
 d^1 = c \oplus^{\theta}_{B^1 \setminus B^2} b^1
 \quad \mbox{and} \quad d^2 = d^1 \oplus^{\theta}_{B^2 \setminus B^1}
 b^2
 $$
 are in $\ang{\cal Q}$, and
 $$
 d_2 = a \oplus^{\theta}_A (\lambda b^1 + (1 - \lambda) b^2) \, .
 $$
 Together with property (iii) of Remark \ref{remprop}.2, this
 implies that
 $$
 a \oplus^{\theta}_A b \in \ang{\cal Q}
 $$
 for all $a,b \in \ang{\cal Q}$, every finite $({\cal F}_{t})$-stopping time
 $\theta$ and $A \in {\cal F}_{\theta}$.
 \end{proof}

 \subsection{Time-consistent coherent utility functional processes}

 \begin{theorem} \label{thmcohnecess}
 Let $T \in{\mathbb N} \cup \{\infty\}$ and
 $(\phi_{t,T})_{t \in [0,T] \cap \mathbb{N}}$
 a relevant time-consistent coherent utility process such that
 $\phi_{0,T}$ can be represented as
 $$
 \phi_{0,T}(X) = \inf_{a \in {\cal Q}} \ang{X,a}_{0,T} \, , \quad
 X \in {\cal R}^{\infty}_{0,T} \, ,
 $$
 for some $\sigma({\cal R}^{\infty} , {\cal A}^1)$-closed,
 convex subset ${\cal Q}$ of ${\cal D}_{0,T}$.
 Then,\\
 {\bf 1.} For every finite $({\cal F}_{t})$-stopping time $\tau \le T$,
 $$
 \phi_{\tau,T}(X) = {\rm ess\,inf }_{a \in {\cal Q}}
 \frac{\ang{X,a}_{\tau,T}}{\ang{1,a}_{\tau,T}} =
 {\rm ess\,inf }_{a \in {\cal Q}^e}
 \frac{\ang{X,a}_{\tau,T}}{\ang{1,a}_{\tau,T}} \, , \quad
 X \in {\cal R}^{\infty}_{\tau,T} \, ,
 $$
 where
 $$
 \frac{\ang{X,a}_{\tau,T}}{\ang{1,a}_{\tau,T}}
 \quad \mbox{is understood to be } \infty \mbox{ on }
 \crl{\ang{1,a}_{\tau,T} = 0} \, ,
 $$
 and
 $$
 {\cal Q}^{e} := {\cal Q} \cap {\cal D}_{0,T}^e \, .
 $$
 {\bf 2.} ${\cal Q}$ and ${\cal Q}^e$ are stable under concatenation.
 \end{theorem}
 \begin{proof} $\mbox{}$

 1. Let $({\cal C}_{t,T})_{t \in [0,T] \cap \mathbb{N}}$ be the
 acceptance set process corresponding to
 $(\phi_{t,T})_{t \in [0,T] \cap \mathbb{N}}$.
 Parts 1 and 2 of Proposition \ref{cthetactau} imply that
 for every finite $({\cal F}_t)$-stopping time $\tau \le T$ and
 $X \in {\cal R}^{\infty}_{\tau,T}$,
 $$
 X \in {\cal C}_{\tau,T} \; \Leftrightarrow \;
 1_A X \in {\cal C}_{0,T} \quad \mbox{for all } A \in {\cal
 F}_{\tau} \, .
 $$
 It follows from Corollary \ref{cohrep} that
 $$
 {\cal Q} = \crl{a \in {\cal D}_{0,T} \mid \phi^{\#}_{0,T}(a) = 0} \, ,
 $$
 and from Corollary \ref{cohrepde} that
 $$
 \phi_{0,T}(X) = \inf_{a \in {\cal Q}^e} \ang{X,a}_{0,T} \,
 ,
 $$
 where ${\cal Q}^e = {\cal Q} \cap {\cal D}^e_{0,T}$.
 Hence, for all $X \in {\cal R}^{\infty}_{\tau, T}$,
 \beas
 X \in {\cal C}_{\tau,T} \; &\Leftrightarrow& \;
 \ang{1_A X, a}_{0,T} \ge 0
 \quad \mbox{for all } A \in {\cal F}_{\tau} \mbox{ and }
 a \in {\cal Q}^e\\
 &\Leftrightarrow&
 \ang{X,a}_{\tau,T} \ge 0 \quad \mbox{for all } a \in {\cal Q}^e
 \, .
 \eeas
 This shows that $\phi_{\tau,T}$ and the coherent utility functional
 $$
 \mbox{ess\,inf }_{a \in {\cal Q}^e}
 \frac{\ang{X,a}_{\tau,T}}{\ang{1,a}_{\tau,T}}
 \, , \quad X \in {\cal R}^{\infty}_{\tau,T}
 $$
 have the same acceptance set. Hence, they must be equal. It is
 clear that
 $$
 {\rm ess\,inf }_{a \in {\cal Q}}
 \frac{\ang{X,a}_{\tau,T}}{\ang{1,a}_{\tau,T}} \le
 {\rm ess\,inf }_{a \in {\cal Q}^e}
 \frac{\ang{X,a}_{\tau,T}}{\ang{1,a}_{\tau,T}} \, , \quad \mbox{for all }
 X \in {\cal R}^{\infty}_{\tau,T} \, .
 $$
 On the other hand, since $X - \phi_{\tau,T}(X) 1_{[\tau, \infty)} \in {\cal
 C}_{\tau,T}$,
 it follows that
 $$
 \ang{1_A \brak{X - \phi_{\tau,T}(X) 1_{[\tau, \infty)}},a}_{0,T} \ge 0 \, , \quad
 \mbox{for all } A \in {\cal F}_{\tau} \mbox{ and } a \in {\cal
 Q} \, ,
 $$
 and therefore,
 $$
 \frac{\ang{\brak{X - \phi_{\tau,T}(X) 1_{[\tau, \infty)}},a}_{\tau,T}}{\ang{1,a}_{\tau,T}}
 \ge 0 \, , \quad \mbox{for all } a \in {\cal Q} \, ,
 $$
 which shows that
 $$
 {\rm ess\,inf }_{a \in {\cal Q}}
 \frac{\ang{X,a}_{\tau,T}}{\ang{1,a}_{\tau,T}}
 \ge \phi_{\tau,T} (X) \, , \quad \mbox{for all }
 X \in {\cal R}^{\infty}_{\tau, T} \, .
 $$

 2. To show that ${\cal Q}$ is stable under concatenation,
 we assume by way of contradiction that there exist $a,b \in {\cal Q}$,
 a finite $({\cal F}_{t})$-stopping time $\theta\le T$ and $A \in {\cal F}_{\theta}$
 such that
 $c := a \oplus^{\theta}_A b \notin {\cal Q}$. Since
 ${\cal Q}$ is $\sigma({\cal R}^{\infty}, {\cal
 A}^1)$-closed and convex, it follows from the
 separating hyperplane theorem that
 there exists an $X \in {\cal R}^{\infty}_{0,T}$ such that
 $$
 \ang{X,c}_{0,T} < \inf_{d \in {\cal Q}}
 \ang{X,d}_{0,T} = \phi_{0,T}(X)
 \, .
 $$
 Note that
 \beas
 &&
 \E{\sum_{j \in [\theta,T] \cap \mathbb{N}} X_j \Delta c_j}\\
 &=& \E{
 1_{\crl{A^c \cup \crl{\ang{1,b}_{\theta,T} = 0}}}
 \sum_{j\in [\theta,T] \cap \mathbb{N}} X_j \Delta a_j
 + 1_{\crl{A \cap \crl{\ang{1,b}_{\theta,T} > 0}}}
 \frac{\ang{1,a}_{\theta,T}}{\ang{1,b}_{\theta,T}}
 \sum_{j\in [\theta,T] \cap \mathbb{N}} X_j \Delta b_j}\\
 &=&
 E\Bigg[\bigg(1_{\crl{A^c \cup \crl{\ang{1,b}_{\theta,T} = 0}}}
 \frac{\ang{X,a}_{\theta,T}}{\ang{1,a}_{\theta,T}} 1_{\crl{\ang{1,a}_{\theta,T} >
 0}}\\
 && \qquad\qquad\qquad\qquad + 1_{\crl{A \cap \crl{\ang{1,b}_{\theta,T} > 0}}}
 \frac{\ang{X,b}_{\theta,T}}{\ang{1,b}_{\theta,T}}
 \bigg) \sum_{j \in [\theta,T] \cap \mathbb{N}} \Delta a_j\Bigg]\\
 &\ge& \E{\phi_{\theta,T}(X) \sum_{j\in
 [\theta,T] \cap \mathbb{N}} \Delta a_j}
 \eeas
 Hence,
 \beas
 \phi_{0,T}(X) &>& \ang{X,c}_{0,T}
 = \E{\sum_{j\in [0,T] \cap \mathbb{N}} X_j \Delta c_j}
 \ge \E{\sum_{j \in [0,\theta)} X_j \Delta a_j + \sum_{j \in
 [\theta,T] \cap \mathbb{N}} \phi_{\theta,T}(X)
 \Delta a_j}\\
 &\ge& \phi_{0,T} \brak{X 1_{[0, \theta)} + \phi_{\theta,T}(X) 1_{[\theta,\infty)}} \,
 ,
 \eeas
 and therefore, $(\phi_{t,T})_{t \in [0,T] \cap \mathbb{N}}$ is
 not time-consistent. But this
 contradicts the assumptions, and therefore, ${\cal Q}$
 has to be stable under concatenation. It follows immediately that
 ${\cal Q}^e$ is stable under concatenation.
 \end{proof}

 \begin{Remark} \label{remcohreptau}
 Let $T \in {\mathbb N}\cup\{\infty\}$ and
 ${\cal Q}^e$ a non-empty subset of ${\cal D}^e_{0,T}$.
 Define for all $t \in [0,T] \cap \mathbb{N}$,
 $$
 \phi_{t,T}(X) := {\rm ess\,inf }_{a \in {\cal Q}^e}
 \frac{\ang{X,a}_{t,T}}{\ang{1,a}_{t,T}} \, , \quad X \in {\cal R}^{\infty}_{t,T} \, .
 $$
 Then, obviously,
 $(\phi_{t,T})_{t \in [0,T] \cap \mathbb{N}}$ is a relevant
 coherent utility functional process, and it is easy to see
 that for every finite $({\cal F}_t)$-stopping time $\tau \le T$
 and all $X \in {\cal R}^{\infty}_{\tau,T}$,
 $$
 \phi_{\tau,T}(X) = {\rm ess\,inf }_{a \in {\cal Q}^e}
 \frac{\ang{X,a}_{\tau,T}}{\ang{1,a}_{\tau,T}} \, .
 $$
 \end{Remark}

 \begin{theorem} \label{thmcohsuff}
 Let $T \in {\mathbb N}\cup\{\infty\}$ and
 ${\cal Q}^e$ a non-empty subset of ${\cal D}^e_{0,T}$
 that is stable under concatenation.
 Define for all $t \in [0,T] \cap \mathbb{N}$,
 $$
 \phi_{t,T}(X) := {\rm ess\,inf }_{a \in {\cal Q}^e}
 \frac{\ang{X,a}_{t,T}}{\ang{1,a}_{t,T}} \, , \quad X \in {\cal R}^{\infty}_{t,T} \, .
 $$
 Then $(\phi_{t,T})_{t \in [0,T] \cap \mathbb{N}}$ is a relevant
 time-consistent coherent utility functional process.
 \end{theorem}
 \begin{proof}
 By Remark \ref{remcohreptau}, $(\phi_{t,T})_{t \in [0,T] \cap
 \mathbb{N}}$ is a relevant coherent utility functional process such that
 for every finite $({\cal F}_t)$-stopping time $\tau \le T$ and
 all $X \in {\cal R}^{\infty}_{\tau,T}$,
 $$
 \phi_{\tau,T}(X) = {\rm ess\,inf }_{a \in {\cal Q}^e}
 \frac{\ang{X,a}_{\tau,T}}{\ang{1,a}_{\tau,T}} \, .
 $$
 To show time-consistency, we denote by $({\cal C}_{t,T})_{t \in
 [0,T] \cap \mathbb{N}}$ the acceptance set process corresponding
 to $(\phi_{t,T})_{t \in [0,T] \cap \mathbb{N}}$ and prove that
 $$
 {\cal C}_{t,T} = {\cal C}_{t, \theta} + {\cal C}_{\theta,T} \, ,
 $$
 for all $t \in [0,T] \cap \mathbb{N}$ and every finite
 $({\cal F}_t)$-stopping time $\theta$ such that $t \le \theta \le
 T$.
 If $Y \in {\cal C}_{t,\theta} \subset {\cal C}_{t,T}$
 and $Z \in {\cal C}_{\theta,T}$, then
 $$
 \ang{Z,a}_{\theta,T} \ge 0 \quad \mbox{for all } a \in {\cal Q}^e \, ,
 $$
 which implies
 $$
 \ang{Z,a}_{t,T} \ge 0 \quad \mbox{for all } a \in {\cal Q}^e \, .
 $$
 Therefore, $Z \in {\cal C}_{t,T}$ and $Y + Z \in {\cal
 C}_{t,T}$. This shows that ${\cal C}_{t,\theta} + {\cal C}_{\theta,T}
 \subset {\cal C}_{t,T}$.
 To prove ${\cal C}_{t,T} \subset {\cal C}_{t,\theta} + {\cal C}_{\theta,T}$,
 we choose a process $X \in {\cal C}_{t,T}$. Obviously, $Z := (X - \phi_{\theta,T}(X))
 1 _{[\theta,\infty)} \in {\cal C}_{\theta,T}$, and it remains to show that
 \beq \label{yinc}
 Y := X - Z = X 1_{[t,\theta)} + \phi_{\theta,T}(X) 1_{[\theta,\infty)}
 \in {\cal C}_{\tau,T} \, .
 \eeq
 Since ${\cal Q}^e$ is stable under concatenation, the set
 $$
 \crl{\frac{\ang{X,a}_{\theta,T}}{\ang{1,a}_{\theta,T}} \mid a \in {\cal Q}^e}
 $$
 is directed downwards. Therefore,
 there exists a sequence $(b^n)_{n \in \mathbb{N}}$ in
 ${\cal Q}^e$ such that
 $$
 \frac{\ang{X,b^0}_{\theta,T}}{\ang{1,b^0}_{\theta,T}} \le
 \N{X}_{\theta, \infty} \le \N{X}_{\infty}
 $$
 and
 $$
 \frac{\ang{X,b^n}_{\theta,T}}{\ang{1,b^n}_{\theta,T}} \searrow
 \phi_{\theta,T}(X) \quad \mbox{almost surely} \, .
 $$
 Moreover,
 $a \oplus^{\theta}_{\Omega} b^n \in {\cal Q}^e$
 for all $a \in {\cal Q}^e$ and $n \in \mathbb{N}$.
 Hence,
 $$
 0 \le \ang{X,a \oplus^{\theta}_{\Omega} b^n}_{t,T}
 \searrow \ang{Y,a}_{t,T} \quad \mbox{almost surely} \, ,
 $$
 and therefore,
 $$
 \ang{Y,a}_{t,T} \ge 0 \quad \mbox{for all }
 a \in {\cal Q}^e \, ,
 $$
 which completes the proof.
 \end{proof}

 \subsection{Time-consistent concave monetary utility functional processes}

 \begin{definition}
 Let $f \in L^0({\cal F})$ and $\tau$ a finite
 $({\cal F}_t)$-stopping time. If there exists a $g \in L^1({\cal
 F})$ such that $f \ge g$, we define
 $$
 \E{f \mid {\cal F}_{\tau}}
 := \lim_{n \to \infty} \E{f \wedge n \mid {\cal F}_{\tau}}
 \, .
 $$
 If there exists a $g \in L^1({\cal F})$ such that $f \le g$, we define
 $$
 \E{f \mid {\cal F}_{\tau}}
 := \lim_{n \to - \infty} \E{f \vee n \mid {\cal F}_{\tau}}
 \, .
 $$
 If $X$ is an adapted process on $(\Omega, {\cal F},
 ({\cal F}_t)_{t \in \mathbb{N}}, P)$ taking values in the
 interval $[m,\infty]$ for some $m \in \mathbb{R}$, we define
 for all $a \in {\cal A}^1_+$,
 $$
 \ang{X,a}_{\tau,\theta}
 := \lim_{n \to \infty} \ang{X \wedge n,a}_{\tau,\theta}
 \, .
 $$
 If $X$ takes values in $[- \infty, m]$, we define
 for all $a \in {\cal A}^1_+$,
 $$
 \ang{X,a}_{\tau,\theta}
 := \lim_{n \to - \infty} \ang{X \vee n,a}_{\tau,\theta}
 \, .
 $$
 \end{definition}

 \begin{Remark} \label{remconcreptau}
 Let $T \in{\mathbb N} \cup \{\infty\}$
 and $(\phi_{t,T})_{t \in [0,T] \cap \mathbb{N}}$ a
 concave monetary utility functional process such that for
 each $t \in [0,T] \cap \mathbb{N}$, $\phi_{t,T}$ is given by
 $$
 \phi_{t,T}(X) = {\rm ess\,inf}_{a \in {\cal D}_{t,T}}
 \crl{ \ang{X,a}_{t,T} - \gamma_{t,T}(a)} \, , \quad
 X \in {\cal R}^{\infty}_{t,T} \, ,
 $$
 for a a special penalty function $\gamma_{t,T}$
 on ${\cal D}_{t,T}$. Then it can easily be checked that
 for all finite $({\cal F}_t)$-stopping times $\tau \le T$,
 $$
 \phi_{\tau,T}(X) = {\rm ess\,inf}_{a \in {\cal D}_{\tau,T}}
 \crl{ \ang{X,a}_{\tau,T} - \gamma_{\tau,T}(a)} \, ,
 \quad X \in {\cal R}^{\infty}_{\tau,T} \, ,
 $$
 where $\gamma_{\tau,T}$ is the special penalty function on
 ${\cal D}_{\tau,T}$ given by
 \beq \label{gammatau}
 \gamma_{\tau,T}(a) := \sum_{t \in [0,T] \cap \mathbb{N}}
 1_{\crl{\tau = t}} \gamma_t(1_{\crl{\tau = t}} a +
 1_{\crl{\tau \neq t}} 1_{[t,\infty)}) \, , \quad
 a \in {\cal D}_{\tau,T} \, .
 \eeq
 \end{Remark}

 \begin{theorem} \label{thmconcnecess}
 Let $T \in \mathbb{N}  \cup \{\infty\}$ and
 $(\phi_{t,T})_{t \in [0,T] \cap \mathbb{N}}$ a time-consistent
 concave monetary utility process such that for all
 $t \in [0,T] \cap \mathbb{N}$,
 $$
 \phi_{t,T}(X) = {\rm ess\,inf}_{a \in {\cal D}_{t,T}}
 \crl{\ang{X,a}_{t,T} - \phi^{\#}_{t,T}(a)} \, , \quad
 X \in {\cal R}^{\infty}_{t,T} \, .
 $$
 Then
 \beq \label{concnecess}
 \phi^{\#}_{\tau,T}(a) =
 {\rm ess\,sup}_{b\in {\cal D}_{\theta,T}}
 \phi^{\#}_{\tau,T} \brak{a \oplus^{\theta}_{\Omega} b}
 + \E{\phi^{\#}_{\theta,T}(a) \mid {\cal F}_{\tau}}
 \, ,
 \eeq
 for every pair of finite $({\cal F}_{t})$-stopping times
 $\tau, \theta$ such that $0 \le \tau \le \theta \le T$ and
 all $a \in {\cal D}_{\tau, T}$.
 \end{theorem}

 \begin{proof}
 Let $\tau$ and $\theta$ be two finite $({\cal F}_t)$-stopping times such that
 $0 \le \tau \le \theta \le T$, and $({\cal C}_{t,T})_{t \in [0,T] \cap \mathbb{N}}$
 the acceptance set process
 corresponding to $(\phi_{t,T})_{t \in [0,T] \cap \mathbb{N}}$.
 It follows from Remark \ref{remtimeconstau}.2 and
 Proposition \ref{decomp} that for all $a \in {\cal D}_{\tau,T}$,
 \bea \label{deca}
 \nonumber
 \phi_{\tau,T}^{\#}(a) &=&
 {\rm ess\,inf}_{X \in {\cal C}_{\tau,T}} \ang{X,a}_{\tau,T}\\
 \nonumber
 &=& {\rm ess\,inf}_{X \in {\cal C}_{\tau,\theta}} \ang{X,a}_{\tau,T}
 + {\rm ess\,inf}_{X \in {\cal C}_{\theta,T}} \ang{X,a}_{\tau,T}\\
 \nonumber
 &=& {\rm ess\,inf}_{X \in {\cal C}_{\tau,\theta}} \ang{X,a}_{\tau,T}
 + \E{{\rm ess\,inf}_{X \in {\cal C}_{\theta,T}}
 \ang{X,a}_{\theta,T} \mid {\cal F}_{\tau}}\\
 &=& {\rm ess\,inf}_{X \in {\cal C}_{\tau,\theta}} \ang{X,a}_{\tau,T}
 + \E{\phi^{\#}_{\theta,T}(a) \mid {\cal F}_{\tau}} \, .
 \eea
 and for all $a \in {\cal D}_{\tau,T}$ and $b \in {\cal D}_{\theta,T}$,
 \beas
 \phi^{\#}_{\tau,T}(a \oplus^{\theta}_{\Omega} b) &=&
 {\rm ess\,inf}_{X \in {\cal C}_{\tau,T}} \ang{X,a \oplus^{\theta}_{\Omega}
 b}_{\tau,T}\\
 &=& {\rm ess\,inf}_{X \in {\cal C}_{\tau,\theta}} \ang{X,a \oplus^{\theta}_{\Omega} b}_{\tau,T}
 + {\rm ess\,inf}_{X \in {\cal C}_{\theta,T}} \ang{X,a \oplus^{\theta}_{\Omega} b}_{\tau,T}\\
 &=& {\rm ess\,inf}_{X \in {\cal C}_{\tau,\theta}} \ang{X,a}_{\tau,T}
 + \E{{\rm ess\,inf}_{X \in {\cal C}_{\theta,T}} \ang{X,b}_{\theta,T}
 \ang{1,a}_{\theta,T} \mid {\cal F}_{\tau}}\\
 &=& {\rm ess\,inf}_{X \in {\cal C}_{\tau,\theta}} \ang{X,a}_{\tau,T}
 + \E{\phi^{\#}_{\theta,T}(b) \ang{1,a}_{\theta,T} \mid {\cal F}_{\tau}} \, .
 \eeas
 By Remark \ref{remconcreptau},
 $$
 \phi_{\theta,T}(X) = {\rm ess\,inf}_{a \in {\cal D}_{\theta,T}}
 \crl{\ang{X,a}_{\theta,T} - \phi^{\#}_{\theta,T}(a)} \, , \quad
 X \in {\cal R}^{\infty}_{\theta,T} \, ,
 $$
 which implies
 $$
 {\rm ess\,sup}_{b \in {\cal D}_{\theta,T}} \phi^{\#}_{\theta,T}(b) = 0
 \, ,
 $$
 and therefore,
 $$
 {\rm ess\,sup}_{b \in {\cal D}_{\theta,T}}
 \E{\phi^{\#}_{\theta,T}(b) \ang{1,a}_{\theta,T} \mid {\cal F}_{\tau}} = 0
 \, ,
 $$
 Hence,
 $$
 {\rm ess\,sup}_{b\in {\cal D}_{\theta,T}}
 \phi^{\#}_{\tau,T} \brak{a \oplus^{\theta}_{\Omega} b}
 = {\rm ess\,inf}_{X \in {\cal C}_{\tau,\theta}} \ang{X,a}_{\tau,T}
 \, ,
 $$
 which together with \eqref{deca}, proves \eqref{concnecess}.
 \end{proof}

 \begin{corollary} \label{corconcnecess}
 Let $T \in \mathbb{N}  \cup \{\infty\}$ and
 $(\phi_{t,T})_{t \in [0,T] \cap \mathbb{N}}$ a relevant time-consistent
 concave monetary utility process such that
 $\phi_{0,T}$ is continuous for bounded decreasing sequences.
 Then
 $$
 \phi_{\tau,T}(X) = {\rm ess\,inf}_{a \in {\cal D}_{\tau,T}}
 \crl{\ang{X,a}_{\tau,T} - \phi^{\#}_{\tau,T}(a)}
 = {\rm ess\,inf}_{a \in {\cal D}_{\tau,T}^e}
 \crl{\ang{X,a}_{\tau,T} - \phi^{\#}_{\tau,T}(a)} \, ,
 $$
 for every finite $({\cal F}_{t})$-stopping time $\tau \le T$, and
 \beas
 \phi^{\#}_{\tau,T}(a) &=&
 {\rm ess\,sup}_{b\in {\cal D}_{\theta,T}}
 \phi^{\#}_{\tau,T} \brak{a \oplus^{\theta}_{\Omega} b}
 + \E{\phi^{\#}_{\theta,T}(a) \mid {\cal F}_{\tau}}\\
 &=& {\rm ess\,sup}_{b \in {\cal D}_{\theta,T}^{e}}
 \phi^{\#}_{\tau,T} \brak{a \oplus^{\theta}_{\Omega} b}
 +  \E{\phi^{\#}_{\theta,T}(a) \mid {\cal F}_{\tau}}
 \, ,
 \eeas
 for every pair of finite
 $({\cal F}_{t})$-stopping times $\tau, \theta$ such that
 $0 \le \tau \le \theta \le T$ and all $a \in {\cal D}_{\tau,
 \theta}$.
 \end{corollary}
 \begin{proof}
 By Theorem \ref{thmrep}, ${\cal C}_{0,T}$ is
 $\sigma({\cal R}^{\infty},{\cal A}^{1})$-closed.
 Let $\tau \le T$ be a finite $({\cal F}_{t})$-stopping time and
 $(X^{\mu})_{\mu \in M}$ a net in ${\cal C}_{\tau,T}$ such
 that $X^{\mu} \to X$ in $\sigma({\cal R}^{\infty},{\cal A}^{1})$ for
 some $X \in {\cal R}^{\infty}_{\tau,T}$. Then, for each $A \in{\cal F}_{\tau}$,
 $1_{A} X^{\mu} \to 1_{A} X$ in $\sigma({\cal R}^{\infty},{\cal A}^{1})$,
 and by Proposition \ref{cthetactau}.1,
 $(1_A X^{\mu})_{\mu \in M}$ is a net in ${\cal C}_{0,T}$. Hence,
 $1_A X \in {\cal C}_{0,T}$, which by Proposition \ref{cthetactau}.2, implies that
 $X \in {\cal C}_{\tau, T}$. This shows that ${\cal C}_{\tau, T}$
 is $\sigma({\cal R}^{\infty}, {\cal A}^1)$-closed. Hence, it
 follows from Theorem \ref{thmrep} that
 \beq \label{normrep}
 \phi_{\tau,T}(X) = {\rm ess\,inf}_{a \in {\cal D}_{\tau,T}}
 \crl{\ang{X,a}_{\tau,T} - \phi^{\#}_{\tau,T}(a)} \, .
 \eeq
 By Proposition \ref{cthetactau}.3, $\phi_{\tau,T}$ is
 $T$-relevant, which by Corollary \ref{correpde}, implies
 that
 \beq \label{erep}
 \phi_{\tau,T}(X) = {\rm ess\,inf}_{a \in {\cal D}_{\tau,T}^e}
 \crl{\ang{X,a}_{\tau,T} - \phi^{\#}_{\tau,T}(a)} \, .
 \eeq
 By Theorem \ref{thmconcnecess}, it follows from \eqref{normrep}
 that
 $$
 \phi^{\#}_{\tau,T}(a) =
 {\rm ess\,sup}_{b\in {\cal D}_{\theta,T}}
 \phi^{\#}_{\tau,T} \brak{a \oplus^{\theta}_{\Omega} b}
 + \E{\phi^{\#}_{\theta,T}(a) \mid {\cal F}_{\tau}}
 \, ,
 $$
 for every pair of finite
 $({\cal F}_{t})$-stopping times $\tau, \theta$ such that
 $0 \le \tau \le \theta \le T$ and all $a \in {\cal D}_{\tau,
 \theta}$.
 In the proof of Theorem \ref{thmconcnecess} we showed that for all
 $a \in {\cal D}_{\tau,T}$ and $b \in {\cal D}_{\theta,T}$,
 $$
 \phi^{\#}_{\tau,T}(a \oplus^{\theta}_{\Omega} b) =
 {\rm ess\,inf}_{X \in {\cal C}_{\tau,\theta}} \ang{X,a}_{\tau,T}
 + \E{\phi^{\#}_{\theta,T}(b) \ang{1,a}_{\theta,T} \mid {\cal F}_{\tau}} \,
 ,
 $$
 and it follows from \eqref{normrep} and \eqref{erep} that
 $$
 {\rm ess\,sup}_{b \in {\cal D}_{\theta,T}} \phi^{\#}_{\theta,T}(b)
 = {\rm ess\,sup}_{b \in {\cal D}_{\theta,T}^e} \phi^{\#}_{\theta,T}(b) = 0
 \, .
 $$
 Hence,
 $$
 {\rm ess\,sup}_{b \in {\cal D}_{\theta,T}}
 \E{\phi^{\#}_{\theta,T}(b) \ang{1,a}_{\theta,T} \mid {\cal F}_{\tau}}
 = {\rm ess\,sup}_{b \in {\cal D}_{\theta,T}^e}
 \E{\phi^{\#}_{\theta,T}(b) \ang{1,a}_{\theta,T} \mid {\cal F}_{\tau}} = 0
 \, ,
 $$
 and therefore,
 $$
 {\rm ess\,sup}_{b\in {\cal D}_{\theta,T}}
 \phi^{\#}_{\tau,T} \brak{a \oplus^{\theta}_{\Omega} b}
 = {\rm ess\,sup}_{b \in {\cal D}_{\theta,T}^{e}}
 \phi^{\#}_{\tau,T} \brak{a \oplus^{\theta}_{\Omega} b} \, .
 $$
 \end{proof}

 \begin{definition}
 Let $T \in \mathbb{N} \cup \infty$ and $\theta \le T$ a finite
 $({\cal F}_t)$-stopping time.\\[2mm]
 For every $a \in {\cal D}_{0,T}$, we define the process
 $\stackrel{\rightarrow \theta}{a} \in {\cal D}_{\theta,T}$ as follows:
 $$
 \stackrel{\rightarrow \theta}{a} :=
 \left\{
 \begin{array}{ll}
 \frac{a}{\ang{1,a}_{\theta,T}} 1_{[\theta,\infty)}
 & \mbox{ on } \crl{\ang{1,a}_{\theta,T} > 0}\\
 1_{[\theta,\infty)} & \mbox{ on } \crl{\ang{1,a}_{\theta,T} = 0}
 \end{array}
 \right.
 $$
 If $\gamma_{\theta,T}$ is a special penalty function on
 ${\cal D}_{\theta,T}$, we extend it to
 ${\cal D}_{0, T}$ by setting
 $$
 \gamma_{\theta,T}^{\rm ext}(a) :=
 \left\{
 \begin{array}{ll}
 \ang{1,a}_{\theta,T} \gamma_{\theta,T}
 \brak{\stackrel{\rightarrow \theta}{a}}
 & \mbox{ on } \crl{\ang{1,a}_{\theta,T} > 0}\\
 0 & \mbox{ on } \crl{\ang{1,a}_{\theta,T} = 0}
 \end{array}
 \right. \, , \, a \in {\cal D}_{0,T} \, .
 $$
 \end{definition}

 \begin{theorem} \label{thmconcsuff}
 Let $T \in{\mathbb N} \cup \{\infty\}$ and
 $(\phi_{t,T})_{t \in [0,T] \cap \mathbb{N}}$ a concave monetary
 utility process such that for every $t \in [0,T] \cap \mathbb{N}$
 and $X \in {\cal R}^{\infty}_{t,T}$,
 $$
 \phi_{t,T}(X) = {\rm ess\,inf}_{a \in {\cal D}_{t,T}}
 \crl{ \ang{X,a}_{t,T} - \gamma_{t,T}(a)} \, ,
 $$
 for a special penalty function $\gamma_{t,T}$
 on ${\cal D}_{t,T}$. Assume that at least
 one of the following two conditions is satisfied:\\[2mm]
 {\rm (1)} For each $t \in [0,T] \cap \mathbb{N}$ and
 every finite $({\cal F}_t)$-stopping time $\theta$ such that
 $t \le \theta \le T$,
 $$
 \gamma_{t,T}(a) = {\rm ess\,sup}_{b \in {\cal D}_{\theta,T}}
 \gamma_{t,T}(a \oplus^{\theta}_{\Omega} b)
 + \E{\gamma^{\rm ext}_{\theta,T}(a)
 \mid {\cal F}_t} \, , \quad
 \mbox{for all } a \in {\cal D}_{t,T} \, .
 $$
 {\rm (2)} $T \in \mathbb{N}$, and for each $t = 0, \dots, T-1$,
 $$
 \gamma_{t,T}(a) =
 {\rm ess\,sup}_{b \in {\cal D}_{t+1,T}}
 \gamma_{t,T}(a \oplus^{t+1}_{\Omega} b)
 + \E{\gamma^{\rm ext}_{t+1,T}(a) \mid {\cal F}_t} \, ,
 \quad \mbox{for all } a \in {\cal D}_{t,T} \, .
 $$
 Then $(\phi_{t,T})_{t \in [0,T] \cap \mathbb{N}}$ is
 time-consistent.
 \end{theorem}
 \begin{proof}
 Let $t \in [0,T] \cap \mathbb{N}$ and $\theta$ a finite
 $({\cal F}_t)$-stopping time
 such that $0\le t \le \theta\le T$. First, note that
 \beq \label{firstineq}
 \gamma_{t,T}(a) \ge
 {\rm ess\,sup}_{b\in{\cal D}_{\theta,T}}
 \gamma_{t,T}(a \oplus^{\theta}_{\Omega} b)
 + \E{\gamma^{\rm ext}_{\theta,T}(a) \mid {\cal F}_{t}}
 \, , \quad \mbox{for all } a \in {\cal D}_{t,T} \, ,
 \eeq
 implies that for all $a \in {\cal D}_{t,T}$ and $b \in {\cal
 D}_{\theta,T}$,
 $$
 \gamma_{t,T}(a \oplus^{\theta}_{\Omega} b) \ge
 \gamma_{t,T}(a) +
 \E{\gamma_{\theta,T}^{\rm ext}(a \oplus^{\theta}_{\Omega} b) \mid {\cal F}_{t}}
 \, ,
 $$
 and therefore,
 \beas
 && \ang{X_{[t, \theta)} +
 \edg{\ang{X,b}_{\theta,T} - \gamma_{\theta,T}(b)} 1_{[\theta,\infty)},
 a}_{t,T} - \gamma_{t,T}(a)\\
 &=& \ang{X, a \oplus^{\theta}_{\Omega} b}_{t,T}
 - \ang{\gamma_{\theta,T}(b) 1_{[\theta,\infty)} , a}_{t,T} -
 \gamma_{t,T}(a)\\
 &=& \ang{X, a \oplus^{\theta}_{\Omega} b}_{t,T}
 - \E{\gamma_{\theta,T}^{\rm ext}(a \oplus^{\theta}_{\Omega} b) \mid {\cal F}_{t}} -
 \gamma_{t,T}(a)\\
 &\ge& \ang{X, a \oplus^{\theta}_{\Omega} b}_{t,T} -
 \gamma_{t,T}(a \oplus^{\theta}_{\Omega} b) \, ,
 \eeas
 for all $X \in {\cal R}^{\infty}_{t,T}$.
 This shows that
 $$
 \phi_{t,T}(X) \le \phi_{t,T}(X 1_{[t, \theta)} +
 \phi_{\theta,T}(X) 1_{[\theta,\infty]}) \, , \quad
 \mbox{for all } X \in {\cal R}^{\infty}_{t,T} \, .
 $$
 On the other hand, it follows from
 $$
 \gamma_{t,T}(a) \le
 {\rm ess\,sup}_{b\in{\cal D}_{\theta,T}}
 \gamma_{t,T}(a \oplus^{\theta}_{\Omega} b)
 + \E{\gamma^{\rm ext}_{\theta,T}(a) \mid {\cal F}_{t}}
 \, , \quad \mbox{for all }a\in{\cal D}_{t,T}
 $$
 that
 \beas
 && \ang{X,a}_{t,T} - \gamma_{t,T}(a)\\
 &\ge& \ang{X,a}_{t,T} -
 \E{\gamma^{\rm ext}_{\theta,T}(a) \mid {\cal F}_{t}}
 - {\rm ess\,sup}_{b\in{\cal D}_{\theta,T}}
 \gamma_{t,T}(a \oplus^{\theta}_{\Omega} b)\\
 &=& \mbox{ess\,inf}_{b\in{\cal D}_{\theta,T}}\crl{\ang{X 1_{[t, \theta)} +
 \edg{\ang{X, \stackrel{\rightarrow \theta}{a}}_{\theta,T} - \gamma_{\theta,T}
 \brak{\stackrel{\rightarrow \theta}{a}}}
 1_{[\theta,\infty)}, a \oplus^{\theta}_{\Omega} b}_{t,T}-
 \gamma_{t,T}(a \oplus^{\theta}_{\Omega} b)}\\
 &\ge& \mbox{ess\,inf}_{b\in{\cal D}_{\theta,T}}
 \crl{\ang{X 1_{[t, \theta)} +
 \phi_{\theta,T}(X) 1_{[\theta,\infty)},
 a \oplus^{\theta}_{\Omega} b}_{t,T}-\gamma_{t,T}
 (a \oplus^{\theta}_{\Omega} b)}\\
 &\ge& \phi_{t,T}\brak{X 1_{[t,\theta)}+\phi_{\theta,T}(X)
 1_{[\theta,\infty)}} \, ,
 \eeas
 for all $X \in {\cal R}^{\infty}_{t,T}$ and
 $a \in {\cal D}_{t, T}$, which shows that
 $$
 \phi_{t,T}(X) \ge \phi_{t,T}(X 1_{[t, \theta)} +
 \phi_{\theta,T}(X) 1_{[\theta,\infty)}) \, , \quad
 \mbox{for all } X \in {\cal R}^{\infty}_{t,T} \, .
 $$
 This shows that it follows directly from condition (1) that
 $(\phi_{t,T})_{t \in [0,T] \cap \mathbb{N}}$ is time-consistent.
 If condition (2) is satisfied, then $(\phi_{t,T})_{t \in [0,T] \cap \mathbb{N}}$
 is time consistent because it fulfills the assumption
 \eqref{onestep} of Proposition \ref{proponestep}.
 \end{proof}
 Exactly the same arguments as in the proof of Theorem
 \ref{thmconcsuff} yield the following

 \begin{corollary} \label{corconcsuff}
 Let $T \in{\mathbb N} \cup \{\infty\}$ and
 $(\phi_{t,T})_{t \in [0,T] \cap \mathbb{N}}$ a concave monetary
 utility process such that for all $t \in [0,T] \cap \mathbb{N}$
 and $X \in {\cal R}^{\infty}_{t,T}$,
 $$
 \phi_{t,T}(X) = {\rm ess\,inf}_{a \in {\cal D}^e_{t,T}}
 \crl{ \ang{X,a}_{t,T} - \gamma_{t,T}(a)} \, ,
 $$
 for a special penalty function $\gamma_{t,T}$
 on ${\cal D}_{t,T}$. Assume that at least
 one of the following two conditions is satisfied:\\[2mm]
 {\rm (1)}
 For each $t \in [0,T] \cap \mathbb{N}$ and
 every finite $({\cal F}_t)$-stopping time $\theta$ such that
 $t \le \theta \le T$,
 $$
 \gamma_{t,T}(a) =
 {\rm ess\,sup}_{b \in {\cal D}^e_{\theta,T}}
 \gamma_{t,T}(a \oplus^{\theta}_{\Omega} b)
 + \E{\gamma^{\rm ext}_{\theta,T}(a)
 \mid {\cal F}_t} \, , \quad \mbox{for all }
 a \in {\cal D}^e_{t,T} \, .
 $$
 {\rm (2)} $T \in \mathbb{N}$, and for each $t = 0, \dots, T-1$,
 $$
 \gamma_{t,T}(a) =
 {\rm ess\,sup}_{b \in {\cal D}^e_{t+1,T}}
 \gamma_{t,T}(a \oplus^{t+1}_{\Omega} b)
 + \E{\gamma^{\rm ext}_{t+1,T}(a) \mid {\cal F}_t} \, ,
 \quad \mbox{for all } a \in {\cal D}^e_{t,T} \, .
 $$
 Then $(\phi_{t,T})_{t \in [0,T] \cap \mathbb{N}}$ is time-consistent.
 \end{corollary}

 \section{Special cases and examples}

 In much of this section the concept of m-stability plays an important
 role. It appears in Artzner et al. (2002) and Delbaen (2004) and
 under different names also in Engwerda et al. (2002),
 Epstein and Schneider (2003), Wang (2003) and Riedel (2004).

 If $T = \infty$, we denote by ${\cal F}_{\infty}$ the sigma-algebra generated by
 $\bigcup_{t \in \mathbb{N}} {\cal F}_t$.

 \begin{definition}
 For $T \in \mathbb{N} \cup \crl{\infty}$,
 $f,g \in L^1({\cal F}_T)$, a finite $({\cal F}_t)$-stopping
 time $\theta \le T$ and $A \in {\cal F}_{\theta}$, we define
 \beq \label{otimes}
 f \otimes^{\theta}_A g :=
 \left\{
 \begin{array}{ll}
 f & \mbox{on } A^c \cup \crl{\E{g \mid {\cal F}_{\theta}} =0 }\\
 \frac{\E{f \mid {\cal F}_{\theta}}}{\E{g \mid {\cal F}_{\theta}}}
 \, g & \mbox{on } A \cap \crl{\E{g \mid {\cal F}_{\theta}} > 0 }
 \end{array}
 \right. \, ,
 \eeq
 and we call a subset
 ${\cal P}$ of $\crl{f \in L^1({\cal F}_T) \mid
 f \ge 0 \, , \, \E{f} = 1}$ m-stable
 if it contains $f \otimes^{\theta}_A g$ for all $f,g \in {\cal P}$,
 every finite $({\cal F}_t)$-stopping
 time $\theta \le T$ and $A \in {\cal F}_{\theta}$.
 \end{definition}
 Let $T \in \mathbb{N} \cup \crl{\infty}$ and ${\cal P}$ a
 non-empty subset of
 $$
 \crl{f \in L^1({\cal F}_T) \mid f \ge 0 \, , \, \E{f} =1 } \,
 .
 $$
 If for all $s \in [0,T] \cap \mathbb{N}$,
 $A \in {\cal F}_s$ and $f, g \in {\cal P}$,
 $f \otimes^s_A g$ is in the $\sigma(L^1, L^{\infty})$-closed, convex
 hull $\ang{\cal P}$ of ${\cal P}$, then it can be shown as in the proof of
 Proposition \ref{estable} that $\ang{\cal P}$ is m-stable.

 \subsection{Processes of coherent utility functionals that depend on the final value}

 Let $T \in \mathbb{N}$ and ${\cal P}$ a non-empty subset of the set
 $$
 \crl{f \in L^1({\cal F}_T) \mid f \ge 0 \, , \,
 \E{f} =1 } \, .
 $$
 Then
 $$
 {\cal Q}({\cal P}):= \crl{f 1_{[T,\infty)} \mid f \in {\cal P}}
 $$
 is a non-empty subset of ${\cal D}_{0,T}$,
 and the concatenation of two elements
 $$
 a = f 1_{[T,\infty)} \quad \mbox{and} \quad
 b = g 1_{[T,\infty)}
 $$
 in ${\cal Q}({\cal P})$ at an $({\cal F}_{t})$-stopping time
 $\theta \le T$ for a set $A \in {\cal F}_{\theta}$ is equal to
 $$
 \brak{f \otimes^{\theta}_A g} \, 1_{[T,\infty)} \, .
 $$
 This shows that ${\cal Q}({\cal P})$ is stable under concatenation if and
 only if ${\cal P}$ is m-stable.

 If ${\cal P}^e$ is a non-empty subset of
 $$
 \crl{f \in L^1({\cal F}_T) \mid f > 0 \, , \,
 \E{f} =1 } \, ,
 $$
 then ${\cal Q}({\cal P}^e)$ is a non-empty subset of ${\cal
 D}^e_{0,T}$, and
 $$
 \phi_{t,T}(X) := {\rm ess\,inf }_{a \in {\cal Q}({\cal P}^e)}
 \frac{\ang{X,a}_{t,T}}{\ang{1,a}_{t,T}}
 = {\rm ess\,inf }_{f \in {\cal P}^e}
 \frac{\E{f X_T \mid {\cal F}_t}}{\E{f \mid {\cal F}_t}} \, , \,
 t = 0, \dots, T \, , \, X \in {\cal R}^{\infty}_{t,T} \, ,
 $$
 defines a relevant coherent utility functional process.

 If ${\cal P}^e$ is m-stable, it follows from Theorem \ref{thmcohsuff}
 that $(\phi_{t,T})_{t=0}^T$ is time-consistent.
 On the other hand, if $(\phi_{t,T})_{t=0}^T$ is time
 consistent, then by Theorem \ref{thmcohnecess},
 the $\sigma({\cal A}^1, {\cal R}^{\infty})$-closed convex hull of
 ${\cal Q}({\cal P}^e)$ is stable under concatenation, which implies that
 the $\sigma(L^1,L^{\infty})$-closed, convex hull of ${\cal P}^e$
 is m-stable.

 This class of time-consistent coherent utility functional
 processes appears in Artzner et al. (2002), Engwerda et al. (2002) and in a
 continuous-time setup in Delbaen (2004). Rosazza Gianin (2003)
 studies the relation between time-consistent monetary utility functionals
 that depend on real-valued random variables and g-expectations.

 \subsection{Processes of coherent utility functionals defined
 by m-stable sets and worst stopping}
 \label{stoppingtimes}

 Let $T \in \mathbb{N} \cup \{\infty\}$ and
 ${\cal P}^e$ a non-empty m-stable subset of
 $$
 \crl{f \in L^1({\cal F}_T) \mid f > 0 \, , \, \E{f} =1 } \, .
 $$
 For all $t \in [0,T] \cap \mathbb{N}$, define
 $$
 \psi_t(Y) := {\rm ess\,inf}_{f \in {\cal P}^e}
 \frac{\E{f \, Y \mid {\cal F}_t}}{\E{f \mid {\cal F}_t}} \, ,
 \quad Y \in L^{\infty}({\cal F}_T) \, ,
 $$
 and for all $X \in {\cal R}^{\infty}_{t,T}$,
 \beq \label{m-stablephi}
 \phi_{t,T}(X) := {\rm ess\,inf}
 \crl{\psi_t(X_{\xi}) \mid \xi \mbox{ a finite
 $({\cal F}_t)$-stopping time such that $t \le \xi \le T$}} \, .
 \eeq
 Then $(\phi_{t,T})_{t \in [0,T] \cap \mathbb{N}}$ is a
 time-consistent relevant coherent utility functional process.

 To see this, note that $\phi_{0,T}$ is a $T$-relevant coherent
 utility functional on ${\cal R}^{\infty}_{0,T}$ that
 can be represented as
 $$
 \phi_{0,T}(X) = \inf_{a \in {\cal Q}({\cal P}^e)}
 \ang{X,a}_{0,T} \, , \quad X \in {\cal R}^{\infty}_{0,T} \, ,
 $$
 where ${\cal Q}({\cal P}^e)$ is the non-empty subset of
 ${\cal D}_{0,T}$ given by
 $$
 {\cal Q}({\cal P}^e) :=
 \crl{\E{f \mid {\cal F}_{\xi}} 1_{[\xi,\infty)} \mid f \in {\cal
 P}^e \, , \, \xi \le T \mbox{ a finite } ({\cal F}_{t})\mbox{-stopping
 time}} \, .
 $$
 It follows from the Corollaries \ref{cohrep} and \ref{cohrepde}
 that
 $$
 \phi_{0,T}(X) = \inf_{a \in {\cal Q}} \ang{X,a}_{0,T}
 = \inf_{a \in {\cal Q}^e} \ang{X,a}_{0,T}
 \, , \quad X \in {\cal R}^{\infty}_{0,T} \, ,
 $$
 where ${\cal Q}$ is the $\sigma({\cal A}^1, {\cal R}^{\infty})$-closed,
 convex hull of ${\cal Q}({\cal P}^e)$ and ${\cal Q}^e = {\cal Q}
 \cap {\cal D}_{0,T}^e$.

 Let $\theta \le T$ be a finite $({\cal F}_{t})$-stopping time,
 $A \in {\cal F}_{\theta}$ and $a, b$ two processes
 in ${\cal Q}({\cal P}^e)$ of the form
 $$
 a = f_a 1_{[\xi_a, \infty)} \quad \mbox{and} \quad
 b = f_b 1_{[\xi_b, \infty)} \, ,
 $$
 where $\xi_a \le T$ and $\xi_b \le T$ are finite
 $({\cal F}_{t})$-stopping times, $f_a = \E{\hat{f}_a \mid {\cal F}_{\xi_a}}$ and
 $f_b = \E{\hat{f}_b \mid {\cal F}_{\xi_b}}$ for
 $\hat{f}_a, \hat{f}_b \in {\cal P}^e$.
 Then
 \beas
 (a \oplus^{\theta}_A b)_t
 &=& 1_{B^c} \, f_a \, 1_{\crl{t \ge \xi_a}}
 + 1_B \,
 \frac{\E{f_a \mid {\cal F}_{\theta}}}{\E{f_b \mid {\cal
 F}_{\theta}}} \,
 f_b \, 1_{\crl{t \ge \xi_b}}\\
 &=& 1_{B^c} \, f_a \, 1_{\crl{t \ge \xi_a}}
 + 1_B \,
 \frac{\E{\hat{f}_a \mid {\cal F}_{\theta}}}{\E{\hat{f}_b \mid {\cal
 F}_{\theta}}}
 f_b \, 1_{\crl{t \ge \xi_b}}\\
 &=& \E{\hat{f} \mid {\cal F}_{\xi}} \, 1_{\crl{t \ge \xi}} \, ,
 \eeas
 where
 $$
 B = A \cap \crl{t \ge \theta} \cap
 \crl{\xi_b \ge \theta}
 \cap \crl{\xi_a \ge \theta} \in {\cal F}_{\theta \wedge \xi_a \wedge \xi_b} \, ,
 $$
 $$
 \hat{f} = \hat{f}_a \otimes^{\theta}_B \hat{f}_b \quad
 \mbox{and} \quad
 \xi = 1_{B^c} \xi_a + 1_B \xi_b \, .
 $$
 It follows from the m-stability of ${\cal P}^e$
 that ${\cal Q}({\cal P}^e)$ is stable under concatenation.
 Proposition \ref{estable} implies that ${\cal Q}$, and therefore also ${\cal Q}^e$
 are stable under concatenation.
 Hence, it follows from Theorem \ref{thmcohsuff} that
 the sequence of functions $(\tilde{\phi}_{t,T})_{t \in [0,T]
 \cap \mathbb{N}}$ given by
 $$
 \tilde{\phi}_{t,T}(X) := {\rm ess\,inf}_{a \in {\cal Q}^e}
 \frac{\ang{X,a}_{t,T}}{\ang{1,a}_{t,T}} \, ,
 \quad t \in [0,T] \cap \mathbb{N} \, , \,
 X \in {\cal R}^{\infty}_{t,T} \, ,
 $$
 is a time-consistent relevant coherent utility functional process,
 and it can easily be checked that $\tilde{\phi}_{t,T} = \phi_{t,T}$
 for all $t \in [0,T] \cap \mathbb{N}$.

 For finite time horizon $T$, this class of time-consistent coherent
 utility functional processes is also discussed in Artzner et al.
 (2002) and in a continuous-time setup in Delbaen (2004).

 \subsection{Processes of coherent utility functionals that depend on
 the infimum over time}

 Let $T \in \mathbb{N} \cup \crl{\infty}$ and ${\cal P}^e$ a non-empty
 subset of the set
 $$
 \crl{f \in L^1({\cal F}_T) \mid f > 0 \, , \, \E{f} =1 } \, .
 $$
 For all $t \in [0,T] \cap \mathbb{N}$, define
 $$
 \psi_t(Y) := {\rm ess\,inf}_{f \in {\cal P}^e}
 \frac{\E{f \, Y \mid {\cal F}_t}}{\E{f \mid {\cal F}_t}} \, ,
 \quad Y \in L^{\infty}({\cal F}_T) \, ,
 $$
 and
 $$
 \phi_{t,T}(X) := \psi_t \brak{\inf_{s \in [t,T] \cap \mathbb{N}} X_s} \, ,
 \quad X \in {\cal R}^{\infty}_{t,T} \, .
 $$
 Then $(\phi_{t,T})_{t \in [0,T] \cap \mathbb{N}}$ is a
 relevant coherent utility functional process. But even if
 ${\cal P}^e$ is m-stable, $(\phi_{t,T})_{t \in [0,T] \cap \mathbb{N}}$
 is in general not time-consistent.

 For an easy counter-example, consider a probability space of the form
 $\Omega = \crl{\omega_1, \omega_2, \omega_3, \omega_4}$ with
 $P[\omega_j] = \frac{1}{4}$ for all $j=1,\dots,4$. Let $T =2$ and
 assume that the filtration $({\cal F}_t)_{t=0}^2$ is given as
 follows:
 ${\cal F}_0 = \crl{\emptyset, \Omega}$, ${\cal F}_1$ is generated
 by the set $\crl{\omega_1, \omega_2}$ and ${\cal F}_2$ is
 generated by the sets $\crl{\omega_j}$, $j=1 ,\dots,4$.
 If ${\cal P}^e = \crl{1}$. Then, for $t \in \crl{ 0,1,2}$ and $X \in {\cal
 R}^{\infty}_{t,2}$,
 $$
 \phi_{t,2}(X) = \E{\inf_{t \le s \le 2} X_s \mid {\cal F}_t} \, .
 $$
 If $X_0 = 2$, $X_1(\omega_1) = X_1(\omega_2) = 4$, $X_1(\omega_3)
 = X_1(\omega_4) = 1$, $X_2(\omega_1) = 5$, $X_2(\omega_2) = 1$,
 $X_2(\omega_3) = 2$ and $X_2(\omega_4) = -1$, then
 $\phi_{0,2}(X) = \frac{3}{4}$. On the other hand,
 $\phi_{1,2}(X) = \frac{5}{2}$ on $\crl{\omega_1, \omega_2}$ and
 $\phi_{1,2}(X) = 0$ on $\crl{\omega_3, \omega_4}$.
 Hence, $\phi_{0,2}(X 1_{\crl{0}} + \phi_{1,2}(X) 1_{[1,2]}) = 1$.

 \subsection{Processes of monetary risk measures that depend on
 an average over time}

 Let $T \in \mathbb{N} \cup \crl{\infty}$ and ${\cal P}^e$ a non-empty
 subset of the set
 $$
 \crl{f \in L^1({\cal F}_T) \mid f > 0 \, , \, \E{f} = 1} \, .
 $$
 For all $t \in [0,T] \cap \mathbb{N}$, define
 $$
 \psi_t(Y) := {\rm ess\,inf}_{f \in {\cal P}^e}
 \frac{\E{f \, Y \mid {\cal F}_t}}{\E{f \mid {\cal F}_t}} \, ,
 \quad Y \in L^{\infty}({\cal F}_T) \, ,
 $$
 and
 $$
 \phi_{t,T}(X) := \psi_t \brak{\frac{\sum_{s \in [t,T] \cap
 \mathbb{N}} \mu_s X_s}{\sum_{s \in [t,T]
 \cap \mathbb{N}} \mu_s}} \, , \quad X \in {\cal R}^{\infty}_{t,T}
 \, ,
 $$
 where $(\mu_s)_{s \in \mathbb{N}}$ is a sequence of non-negative
 numbers such that
 $$
 \sum_{s \in [0,T] \cap \mathbb{N}} \mu_s = 1 \, ,
 $$
 and
 $$
 \sum_{s \in [t,T] \cap \mathbb{N}} \mu_s > 0 \quad \mbox{for all }
 t \in [0,T] \cap \mathbb{N} \, .
 $$
 Then $(\phi_{t,T})_{t \in [0,T] \cap \mathbb{N}}$ is obviously a
 relevant coherent utility functional process, and if ${\cal P}^e$
 is m-stable, then $(\phi_{t,T})_{t \in [0,T] \cap \mathbb{N}}$
 is time-consistent.

 To see this we denote for $f \in {\cal P}^e$
 by $J(f)$ the process $a \in {\cal D}^e_{0,T}$ given by
 $$
 \Delta a_t := \mu_t \E{f \mid {\cal F}_t} \quad \mbox{for all }
 t \in \mathbb{N} \, .
 $$
 Clearly,
 $$
 \phi_{t,T}(X) = {\rm ess\,inf}_{a \in J({\cal P}^e)}
 \frac{\ang{X,a}_{t,T}}{\ang{1,a}_{t,T}}
 \quad \mbox{for all } t \in [0,T] \cap \mathbb{N} \quad
 \mbox{and} \quad X \in {\cal R}^{\infty}_{t,T} \, ,
 $$
 and it can easily be checked that for all $f,g \in {\cal P}^e$,
 every finite $({\cal F}_t)$-stopping time $\theta \le T$ and
 $A \in {\cal F}_{\theta}$,
 $$
 J(f) \oplus^{\theta}_A J(g) = J(f \otimes^{\theta}_A g) \, .
 $$
 Hence, $J({\cal P}^e)$ is stable under concatenation, and it
 follows from Theorem \ref{thmcohsuff} that
 $(\phi_{t,T})_{t \in [0,T] \cap \mathbb{N}}$ is time-consistent.

 \subsection{Processes of robust entropic utility functionals}

 Let $T \in \mathbb{N}$ and ${\cal P}^e$ a non-empty
 subset of
 $$
 \crl{f \in L^1({\cal F}_T) \mid f>0 \, ; \,
 \E{f} =1} \, .
 $$
 For $t = 0, \dots, T$ and $X \in {\cal R}^{\infty}_{t,T}$, define
 $$
 \phi_{t,T}(X) := \mbox{ess inf}_{f \in{\cal P}^e}
 \crl{- \log \frac{\E{f \exp(- X_T) \mid {\cal F}_t}}{\E{f \mid {\cal
 F}_t}}} \, , \quad X \in {\cal R}^{\infty}_{0,T} \, .
 $$
 Then, for all $t=0, \dots, T$, $\phi_{t,T}$ is a
 $T$-relevant concave monetary utility functional on
 ${\cal R}^{\infty}_{t,T}$ that is
 continuous for bounded decreasing sequences, and
 $(\phi_{t,T})_{t = 0}^T$ is time-consistent if
 ${\cal P}^e$ is m-stable.
 Indeed, it is obvious that for all $t = 0,\dots,T$,
 $\phi_{t,T}$ is a $T$-relevant monetary utility functional on
 ${\cal R}^{\infty}_{t,T}$. To show the other assertions, we introduce for all $f
 \in {\cal P}^e$ and $t = 0, \dots, T$, the mappings
 $$
 \psi^f_t(Y) := - \log \frac{\E{f \exp(-Y) \mid {\cal F}_t}}{\E{f \mid {\cal
 F}_t}} \, , \quad Y \in L^{\infty}({\cal F}_T)
 $$
 and
 $$
 \psi_t(Y) := {\rm ess\,inf} _{f \in {\cal P}^e} \psi^f_t(Y) \, ,
 \quad Y \in L^{\infty}({\cal F}_T) \, .
 $$
 To see that for all $t=0,\dots,T$, $\phi_{t,T}$ is continuous
 for bounded decreasing sequences we let $(Y^n)_{n \in \mathbb{N}}$
 be a decreasing sequence in $L^{\infty}({\cal F}_T)$
 and $Y \in L^{\infty}({\cal F}_T)$
 such that $\lim_{n \to \infty} Y^n = Y$ almost surely.
 Then $\lim_{n \to \infty} \psi_t(Y^n)$ exists almost surely and
 $\lim_{n \to \infty} \psi_t(Y^n) \ge \psi_t(Y)$. On the other
 hand, there exists a sequence $(f^k)_{k \in \mathbb{N}}$ in
 ${\cal P}^e$ such that $\psi_t(Y) = \inf_{k \in \mathbb{N}}
 \psi^{f^k}_t(Y)$ almost surely, and for all $k \in
 \mathbb{N}$,
 $$
 \psi^{f^k}_t(Y) = \lim_{n \to \infty} \psi^{f^k}_t(Y^n)
 \ge \lim_{n \to \infty} \psi_t(Y^n) \, ,
 $$
 Hence, $\psi_t(Y) \ge \lim_{n \to \infty} \psi_t(Y^n)$,
 which shows that $\phi_{t,T}$ is continuous for bounded decreasing
 sequences.

 To see that $\phi_{t,T}$ satisfies condition (3) of Definition
 \ref{def1}, fix an $f \in {\cal P}^e$ and a $Z \in L^{\infty}({\cal F}_T)$.
 Then, define $f_Z \in L^1({\cal F}_T)$ by
 $$
 f_Z := \frac{f \exp \brak{-Z}}{\E{f \exp \brak{-Z} \mid {\cal F}_t}} \, ,
 $$
 and note that it follows from Jensen's inequality that
 for all $Y \in L^{\infty}({\cal F}_T)$,
 \beas
 \psi^f_t(Y) &=&
 - \log \frac{\E{f_Z \frac{f}{f_Z} \exp \brak{-Y} \mid {\cal F}_t}}{\E{
 f \mid {\cal F}_t}}\\
 &=& - \log \E{f_Z \exp \brak{Z-Y} \mid {\cal F}_t}
 - \log \frac{\E{f \exp \brak{-Z} \mid {\cal F}_t}}{\E{
 f \mid {\cal F}_t}}\\
 &\le& \E{f_Z (Y - Z) \mid {\cal F}_t} + \psi^f_t(Z) \, .
 \eeas
 This shows that for all $Y \in L^{\infty}({\cal F}_T)$,
 \beq \label{condentr}
 \psi^f_t(Y) = {\rm ess\,inf}_{
 Z \in L^{\infty}({\cal F}_T)}
 \crl{\E{f_Z Y \mid {\cal F}_t} - \E{f_Z Z \mid {\cal F}_t}
 + \psi^f_t(Z)} \, ,
 \eeq
 and therefore,
 $$
 \psi_t(Y) = {\rm ess\,inf}_{f \in {\cal P}^e \, , \,
 Z \in L^{\infty}({\cal F}_T)}
 \crl{\E{f_Z Y \mid {\cal F}_t} - \E{f_Z Z \mid {\cal F}_t}
 + \psi^f_t(Z)} \, ,
 $$
 from which it can be seen that $\phi_{t,T}$ satisfies condition (3) of Definition
 \ref{def1}.

 Now, assume that ${\cal P}^e$ is m-stable. Then,
 for all $t=0,\dots,T$ and $Y \in L^{\infty}({\cal F}_T)$, the set
 $$
 \crl{\psi^f_t(Y) \mid f \in {\cal P}^e}
 $$
 is directed downwards because for all $f,g \in {\cal P}^e$,
 $$
 \psi^f_t(Y) \wedge \psi^g_t(Y) = \psi^h_t(Y) \, ,
 $$
 where
 $$
 h = f \otimes^t_A g \quad \mbox{for} \quad
 A = \crl{\psi^g_t(Y) < \psi^f_t(Y)} \, .
 $$
 Hence, there exists a sequence $(f^k)_{k
 \in \mathbb{N}}$ in ${\cal P}^e$ such that almost surely,
 $$
 \psi^{f^k}_t(Y) \searrow \psi_t(T) \, , \quad \mbox{as } k
 \to \infty \, .
 $$
 Next, note that for all $t=0,\dots,T-1$, $f,g \in {\cal P}^e$ and $Y \in
 L^{\infty}({\cal F}_T)$,
 $$
 \psi^f_t(\psi^g_{t+1}(Y)) = \psi^h_t(Y) \, ,
 $$
 where
 $$
 h = f \otimes^{t+1}_{\Omega} g \, .
 $$
 It follows that
 $$
 \psi_t(\psi_{t+1}(Y)) = \psi_t(Y) \, ,
 $$
 for all $t=0,\dots,T-1$ and $Y \in L^{\infty}({\cal F}_T)$, and therefore,
 $$
 \phi_{t,T}(X 1_{\crl{t}} + \phi_{t+1,T}(X)
 1_{[t+1,\infty)}) = \phi_{t,T}(X) \, ,
 $$
 for all $t=0,\dots,T-1$ and $X \in {\cal R}^{\infty}_{t,T}$,
 which by Proposition \ref{proponestep}, implies that
 $(\phi_{t,T})_{t=0}^T$ is time-consistent.

 The utility functional $\phi_{0,T}$ is a robust version of the
 mapping
 $$
 C : L^{\infty}({\cal F}_T) \to \mathbb{R} \, , \quad
 Y \mapsto - \log \E{\exp(-Y)} \, ,
 $$
 which assigns a random variable $Y \in L^{\infty}({\cal F}_T)$ its
 certainty equivalent under the exponential utility function
 $$
 x \mapsto - \exp(-x) \, .
 $$
 It is well known that $C$ admits the representation
 \beq \label{entr}
 C(Y) = \inf_{Q} \crl{\mbox{E}_Q[Y] + H(Q \mid P)} \, ,
 \eeq
 where the infimum is taken over all probability measures
 $Q$ on $(\Omega, {\cal F}_T)$ and $H(Q \mid P)$ is the relative entropy of $Q$ with respect
 to $P$. In fact, it can easily be checked that the penalty function in
 the representation \eqref{condentr} is the conditional
 relative entropy of $f_Z$ with respect to $f$.
 For more details and relations to pricing in incomplete markets we refer
 to Frittelli (2000), Rouge et al. (2000) and Delbaen et al. (2002).
 More on the entropic risk measure $-C$ can be found in F\"ollmer and Schied
 (2002a) and Weber (2003).
 A conditional version of the entropic risk measure is studied in Detlefsen (2003).
 In Frittelli and Rosazza Gianin (2004) it is shown that the
 dynamic entropic risk measure is time-consistent in continuous
 time.

 \subsection{Time-consistent monetary utility functional
 processes and worst stopping}

 For finite time horizon the coherent utility functional processes
 of Subsection \ref{stoppingtimes} can be generalized as follows:

 Let $T \in \mathbb{N}$ and $(\hat{\phi}_{t,T})_{t =0}^T$ a
 time-consistent monetary utility functional process on
 $R^{\infty}_{0,T}$ such that for all $t = 0, \dots, T$ and $X \in
 {\cal R}^{\infty}_{t,T}$, $\hat{\phi}_{t,T}(X)$ depends only on
 the final value $X_T$ of $X$, that is, for all $t=0,\dots,T$,
 $$
 \hat{\phi}_{t,T}(X) = \psi_t(X_T) \, ,
 $$
 where $\psi_t$ is a mapping from $L^{\infty}({\cal F}_T)$ to
 $L^{\infty}({\cal F}_t)$ that satisfies the following
 conditions:\\[2mm]
 {\rm (0)} $\psi_t(1_A Y) = 1_A \psi(Y)$ for all
 $Y \in L^{\infty}({\cal F}_T)$ and $A \in {\cal
 F}_t$\\[2mm]
 {\rm (1)} $\psi_t(Y) \le \psi_t(Z)$ for all $Y,Z \in L^{\infty}({\cal F}_T)$
 such that $Y \le Z$\\[2mm]
 {\rm (2)} $\psi_t(Y + m) = \psi_t(Y) + m$ for all $Y \in L^{\infty}({\cal F}_T)$
 and $m \in L^{\infty}({\cal F}_t)$\\[2mm]
 {\rm (tc)} $\psi_t(\psi_{t+1}(Y)) = \psi_t(Y)$, for all
 $t =0 ,\dots, T-1$ and $Y \in L^{\infty}({\cal F}_T)$.\\[2mm]
 Denote by $\Theta_{t,T}$ the set of all
 $({\cal F}_{t})$-stopping times $\xi$ such that
 $t \le \xi \le T$, and define a new monetary utility functional process by
 $$
 \phi_{t,T}(X) :=
 {\rm ess\,inf}_{\xi\in\Theta_{t,T}} \psi_t(X_{\xi}) \,
 .
 $$
 For a given $X \in {\cal R}^{\infty}_{0,T}$, define the process
 $(S_t(X))_{t=0}^T$ recursively by
 $$
 S_T(X) := X_T
 $$
 and
 $$
 S_t(X) := X_t \wedge \psi_t (S_{t+1}(X)) \, ,
 \quad \mbox{for } t \le T-1 \, .
 $$
 For all $t = 0,\dots,T$, denote by $\xi^t$ the stopping time given by
 $$
 \xi^t := \inf \crl{j=t,\dots,T \mid S_j(X) = X_j} \, .
 $$
 It can easily be checked that
 $$
 \phi_{t,T}(X) = \psi_t(X_{\xi^t}) = S_t(X) \, .
 $$
 For a $t \in [0,T-1] \cap \mathbb{N}$, set
 $Y := X 1_{\crl{t}} + \phi_{t,T}(X) 1_{[t+1, \infty)}$. It is easy to see
 that $S_{t+1}(Y) = S_{t+1}(X)$. Therefore,
 $$
 \phi_{t,T}(Y) = S_t(Y) = Y_t \wedge \psi_t (S_{t+1}(Y))
 = X_t \wedge \psi_t (S_{t+1}(X)) = S_t(X) = \phi_{t,T}(X) \, ,
 $$
 which shows that $(\phi_{t,T})_{t=0}^T$ is time-consistent.

 \bigskip \bigskip

 \noindent
 {\bf \Large References}

 \medskip

 \noindent
 Artzner, Ph., Delbaen, F., Eber, J.M., Heath, D. (1997).
 {\sl Thinking coherently}, RISK 10, November, 68--71.\\[2mm]
 Artzner, Ph., Delbaen, F., Eber, J.M., Heath, D.
 (1999).
 {\sl Coherent measures of risk}, Math. Finance {\bf 9}(3), 203--228.\\[2mm]
 Artzner, Ph., Delbaen, F., Eber,
 J.M., Heath, D., Ku, H.
 (2003).
 {\sl Coherent multiperiod risk adjusted values and Bellman's
 principle.} Submitted.\\[2mm]
 Cheridito, P., Delbaen F., Kupper M. (2004a).
 {\sl Coherent and convex risk measures for bounded c\`adl\`ag
 processes}, Stoch. Proc. Appl. {\bf 112}(1), 1--22.\\[2mm]
 Cheridito, P., Delbaen F., Kupper M. (2004b).
 {\sl Coherent and convex risk measures for unbounded c\`adl\`ag
 processes.} Submitted. \\[2mm]
 Cvitani\'{c}, J., Karatzas, J. (1999).
 {\sl On dynamic measures of risk.}
 Finance and Stochastics {\bf 3}(4), 451--482.\\[2mm]
 Delbaen, F. (2002). {\sl Coherent risk measures
 on general probability spaces.}
 Essays in Honour of Dieter Sondermann, Springer-Verlag.\\[2mm]
 Delbaen, F. (2000). Coherent Risk Measures. Cattedra Galileiana. Scuola
 Normale Superiore di Pisa.\\[2mm]
 Delbaen, F. (2003). {\sl The structure of m-stable sets
 and in particular of the set of risk neutral measures.} Working Paper,  ETH Z\"urich.\\[2mm]
 Delbaen, F., Grandits P., Rheinl\"ander T., Samperi D., Schweizer M., Stricker Ch. (2002).
 {\sl Exponential hedging and entropic penalties}, Math. Finance {\bf 12}(2), 99--123.\\[2mm]
 Detlefsen, K. (2003).
 Bedingte und Mehrperiodische Risikomasse. Diplomarbeit,
 Humboldt-Universit\"at zu Berlin.\\[2mm]
 Engwerda J., Roorda B., Schumacher J.M. (2002). {\sl Coherent acceptability measures in
 multiperiod models}, Working Paper, Tilburg University.\\[2mm]
 Eppstein, L.G., Schneider, M. (2003).
 {\sl Recursive multiple-priors}, J. Econom. Theory {\bf 113}(1),
 1--31.\\[2mm]
 F\"ollmer, H., Schied, A. (2002a).
 {\sl Convex measures of risk and trading constraints.} Finance and Stochastics
 {\bf 6}(4), 429--447.\\[2mm]
 F\"ollmer, H., Schied, A. (2002b).
 {\sl Robust preferences and convex measures of risk.} Advances in
 Finance and Stochastics,
 Springer-Verlag.\\[2mm]
 F\"ollmer, H., Schied, A. (2002c).
 Stochastic Finance, An Introduction in Discrete Time.
 de Gruyter Studies in Mathematics 27.\\[2mm]
 Frittelli, M. (2000).
 {\sl Introduction to a theory of value coherent with the no-arbitrage principle.} Finance and
 Stochastics {\bf 4}(3), 275--297 \\[2mm]
 Frittelli, M., Rosazza Gianin, E. (2004).
 {\sl Dynamic convex risk measures.}
 Risk Measures for the 21st Century, Chapter 12,
 Wiley Finance. \\[2mm]
 %\bibitem[FG]{FG} Frittelli, M., Gianin, E.R. (2002).
 %{\sl Putting order in risk measures} Journal of Banking and Finance 26(7), 1473-1486
 Neveu, J. (1975). Discrete-Parameter Martingales.
 North-Holland Publishing Company-Amsterdam, Oxford.\\[2mm]
 Riedel, F. (2004).
 {\sl Dynamic coherent risk measures.} Stoch. Proc. Appl.
 {\bf 112}(2), 185--200.\\[2mm]
 Rosazza Gianin, E. (2003). {\sl Some examples of risk measures via g-expectations},
 Working Paper n. 41, Universit\`a di Milano Bicocca \\[2mm]
 Rouge, R., El Karoui, N.
 (2000).
 {\sl Pricing via utility maximization and entropy}, Mathematical Finance
 {\bf 10}(2), 259--276 \\[2mm]
 Scandolo, G. (2003).
 Risk measures in a dynamic setting. Ph.D. Thesis, Universit\`a degli Studi
 Milano \& Universit\`a di Firenze.\\[2mm]
 Weber, S. (2003).
 {\sl Distribution-invariant dynamic risk measures.} Preprint,
 Humboldt-Universit\"at zu Berlin.\\[2mm]
 Wang, T. (2003).
 {\sl Conditional preferences and updating.} J. Econom Theory {\bf 108}(2),
 286--321.
 \end{document}